\documentclass[11pt]{article}
\usepackage{amsmath,amssymb}
\usepackage{CJK,indentfirst}
\usepackage{amssymb}
\usepackage{amsthm}
\usepackage{cite}
\begin{document}
\baselineskip=20pt
\parskip=1mm
\parindent=20pt
\voffset=-1 true cm \setcounter{page} {1}
\numberwithin{equation}{section}
\newtheorem{Theorem}{\bf Theorem}[section]
\newtheorem{Lem}{\bf Lemma}[section]
\newtheorem{Prop}[Lem]{\bf Proposition}
\newtheorem{Def}{\bf Definition}[section]
\newtheorem{Coro}{\bf Corollary}[section]
\newtheorem{Rem}{\bf Remark}[section]

\title{\textbf{Normalized Solutions to Kirchhoff Equation with Nonnegative Potential}   }
\date{ }
\author{    Shuai Mo\ \ \  Shiwang Ma\thanks{Corresponding author. E-mails: moshuai@126.com(S.Mo)\ \ \  shiwangm@nankai.edu.cn(S.Ma)}  \\
\it \small  School of Mathematical Sciences and LPMC, Nankai University\\
\it \small Tianjin 300071, China }\date{}\maketitle

{\bf \noindent Abstract}\quad This paper is concerned with the existence of solutions to the problem

$$-\left(a+ b\int_{\mathbb{R}^{N}}|\nabla u|^{2} dx  \right)\Delta u +V(x)u+\lambda u = |u|^{p-2}u,\ \  x \in \mathbb{R}^{N},\ \ \lambda \in \mathbb{R}^{+}
                  $$
where $a, b>0$ are constants, $ V \geq 0$ is a potential, $N \geq 1 $, and $ p \in (2+ \frac{4}{N},2^*$). We use a more subtle
analysis to revisit the limited problem($V \equiv 0$), and obtain a new energy inequality and bifurcation results. Based on these observations, we establish the existence of bound state normalized solutions under different assumptions on $V$. These conclusions extend some known results in previous papers.

{\bf \noindent Keywords:} Normalized solutions, Kirchhoff equations, Bound state, Variational methods.

{\bf \noindent 2010 MSC:}\quad 35J50; 35J15; 35J60.

\section{Introduction}

We look for solutions $(u,\lambda) \in H^1(\mathbb{R}^N) \times \mathbb{R}^{+}$ of the problem

$$\left\{\aligned &-\left(a+ b\int_{\mathbb{R}^{N}}|\nabla u|^{2} dx  \right)\Delta u +V(x)u+\lambda u = |u|^{p-2}u, &x \in \mathbb{R}^{N},\\
                  &\int_{\mathbb{R}^{N}} |u|^{2}dx =c^2          \endaligned\right.  \eqno(  K_{V,c}  )$$
where $a, b,c>0$ are constants, $ V: \mathbb{R}^{N}\rightarrow \mathbb{R}^{+} \cup \{0\} $ is a potential, $N \geq 1 $, and $ p \in (2+ \frac{4}{N},2^*$). Moreover, we consider the case $V(x)\rightarrow 0$ as $ |x| \rightarrow \infty $ and allow that $V$ has singularities.

The Kirchhoff problem is related to the stationary analog of the equation
$$  \rho \frac{\partial^2 u}{\partial t^2}-\left(  \frac{P_0}{h}+\frac{E}{2L}\int_{0}^{L} \left|\frac{\partial u}{\partial x}\right|^2 dx   \right) \frac{\partial^2 u}{\partial t^2} =0  $$
presented by Kirchhoff in \cite{Kirchhoff1883} in order to describe the transversal free vibrations of a clamped string in which the dependence of the tension on the deformation cannot be neglected.
We refer readers to \cite{Arosio1997,Bernstein1940,Ca2001,Spagnolo1992} for the physical background on Kirchhoff problem.

Comparing with the corresponding ``local" equations ($b=0$), the presence of the nonlocal term $(\int_{\mathbb{R}^N} |\nabla u|^2 dx )\triangle u $
makes $(K_{V,c})$ is not a pointwise identity, which causes additional mathematical difficulties.
It has received much attention after J.L.Lions proposed an abstract framework in the celebrated paper \cite{J.Lions1978}.  If $\lambda \in \mathbb{R}$ is a fixed parameter, critical points theory is used to look for nontrivial solutions, see \cite{Santos2014,Guo2015,Mao2009,Mao2022,Zou2014}. However, nothing can be given a priori on the $L^2$-norm of the solutions in this method.

There is another way to study problem $(K_{V,c})$, by searching for solutions satisfying $L^2$-norm constraint.
In such a point of view, $u$ with prescribed $L^2$-norm is know as a normalized solution, and $\lambda $ appears as a lagrange multiplier. To obtain such solutions, one usually look for critical points of functional
\begin{equation}\label{equ:functional V}
I(u)=\frac{a}{2}\int_{\mathbb{R}^N}|\nabla u|^2 dx + \frac{1}{2} \int_{\mathbb{R}^N}V(x)u^2 dx + \frac{b}{4} \left(\int_{\mathbb{R}N} |\nabla u|^2 dx\right)^2 - \frac{1}{p}\int_{\mathbb{R}^N} |u|^p dx
\end{equation}
constrained to the $L^2$-sphere
$$ S_c = \left\{  u \in H^1(\mathbb{R}^N) : \int_{\mathbb{R}^N} u^2 dx =c^2   \right\}.  $$
The study on the normalized solutions of partial differential equations has been the purpose of very active research in recent years.
Let us introduce and review some related results in this respect, and classify the influences of potential and nonlocal term.

\noindent\textbf{1.1 Local Case: $b=0$.} For $b=0$, we may assume $a=1$ and $V(x)\rightarrow0$ as $|x|\rightarrow +\infty$. The problem reduces to
$$\left\{\aligned &-\Delta u +V(x)u+\lambda u = |u|^{p-2}u, &x \in \mathbb{R}^{N},\\
                  &\int_{\mathbb{R}^{N}} |u|^{2}dx =c^2.          \endaligned\right.  \eqno(  S_{V,c}  )$$
For  fixed $c>0$, we consider  the minimizing problem
\begin{equation}\label{Minipro}
  e(c):= \inf_{S_c} I(u).
\end{equation}
By the Galiardo-Nirenberg inequality with the best constant \cite{Weinstein},  for  $p \in [2, 2^* )$, we have
\begin{equation}\label{GN ineq}
  \| u \|_p \leq \left( \frac{p}{  2\| Q_p \|_2^{p-2}  }  \right)^{\frac{1}{p}} \| \nabla u\|_2^{\frac{N(p-2)}{2p}}\| u \|_2^{1-\frac{N(p-2)}{2p}},
\end{equation}
with equality only for $u=Q_p$, where up to translations, $Q_p$ is the unique positive solution of
\begin{equation}\label{equ:Q sovles equ}
  -\frac{N(p-2)}{4} \triangle Q_p + \frac{2N-p(N-2)}{4} Q_p=|Q_p|^{p-2}Q_p, \ \ x\in \mathbb{R}^N.
\end{equation}
It is easy to see from (\ref{equ:Q sovles equ}) that
\begin{equation}\label{equ:Q norm}
\| Q_p \|_2^2 = \frac{2}{p}\| Q_p \|_p^p = \|\nabla Q_p \|_2^2.
\end{equation}
In particular,\begin{equation}\label{c GN ineq}
  \| u \|_p \leq \left( \frac{p}{  2\| Q_p \|_2^{p-2}  }  \right)^{\frac{1}{p}} \| \nabla u\|_2^{\frac{N(p-2)}{2p}}c^{1-\frac{N(p-2)}{2p}},\ \ \forall u \in S_c
\end{equation}
with equality only for $u=t \star \left(\frac{c}{\|Q_p\|_2}Q_p\right)$ for all $t\in \mathbb{R}$, where the notation $t \star u(\cdot) = t^{\frac{N}{2}}u(t\cdot)$ for all $u \in H^{1}(\mathbb{R}^N)$ and $t>0$. It is standard to show $e(c) > -\infty$ for any $c>0$ if $2<p<2+\frac{4}{N}$($L^2$-subcritical case), and $e(c) = -\infty$ for any $c>0$ if $2+\frac{4}{N}<p<2^*$($L^2$-supercritical case).
The number  $\bar{p}:= 2+\frac{4}{N}$ is called  $L^2$-critical exponent coming from (\ref{GN ineq}). When dealing with nonlinearities of $L^2$-subcritical case, one could solve (\ref{Minipro}) directly. We refer the readers to \cite{Ikoma,Shibata2014,Zhong2021} and the references therein.

While dealing with the $L^2$-supercritical case involving potential, the functional is unbounded from below, thus the minimizing argument on $S_c$ is not valid anymore. In the classical paper \cite{Jeanjean1997}, L. Jeanjean discussed the problem with $V\equiv 0$
$$\left\{\aligned &-\Delta u + \lambda u = g(u), &x \in \mathbb{R}^{N},\\
                  &\int_{\mathbb{R}^{N}} |u|^{2}dx =c^2,          \endaligned  \right.\eqno{( \tilde{S}_{\infty,c} )} $$
where $g(u)$ is nonhomogeneous $L^2$-supercritical nonlinearity and $G(t)=\int_{0}^{t}g(\tau)d\tau$. By developing a mountain-pass argument on $S_c$, and constructing a $(PS)$ sequence with additional property, L.Jeanjean obtained a
ground state solution $u$ satisfying
$$ \mathcal{I}_{0}(u)=E_{0}(c):=\inf_{v \in \mathcal{P}}\mathcal{I}_{0}(v), $$
where $$\mathcal{I}_{V}(u)=\frac{1}{2}\|\nabla u \|_2^2 + \frac{1}{2}\int_{\mathbb{R}^N}V(x) u^2 dx -\int_{\mathbb{R}^N}G(u)dx $$
and $$ \mathcal{P}=\left\{ u \in S_c : \frac{\partial \mathcal{I}_{V}(t\star u)}{\partial t}\Big|_{t=1}=0 \right\}.$$
The key step is to establish the strict so-called subadditive inequality $E_{0}(a+b)<E_{0}(a)+E_{0}(b)$, which is crucial for using the splitting lemma\cite[Lemma 3.1]{Benci} to recover the compactness. Another method is minimizing $\mathcal{I}_{0}(u)$ over $ \mathcal{P}$, Bartsch, Sovae \cite{BartschJFA} and Yang \cite{yang} proved $\mathcal{I}_{0}(u)$ is coercive on $ \mathcal{P}$, and the minimizing sequence $u_n$ converge strongly to $u\neq0$ in $H^1(\mathbb{R}^N)$ using the strict subadditive inequality. Since $ \mathcal{P}$ is a nature constraint and contains any solutions of $(\tilde{S}_{\infty,c})$, they obtained a ground state solution.

The case of negative potential $V \leq 0$ is considered in \cite{Ding2014} and  \cite{Molle2021}, they both obtained the ground state solution. In fact, under some explicit smallness assumption on $V(x)$, one has $E_{V}(c)<E_{0}(c)$, thus the method in \cite{BartschJFA,Jeanjean1997,yang} is still valid and the trapping nature of the potential provides enough
compactness. However,  when $V\geq 0$, we only have $E_{V}(c)=E_{0}(c)$ and the problem becomes more delicate and difficult.
Although  the mountain pass structure by Jeanjean \cite{Jeanjean1997} is destroyed,   Bartsch et al.\cite{Bartsch2021} established a new variational principle exploiting the Pohozaev identity.
By constructing a suitable linking geometry, the authors in \cite{Bartsch2021} succeeded to  obtain the existence of bound state solutions with high Morse index. But their method only works for power nonlinearity $g(u)=|u|^{p-2}u$ with $L^2$-supercritical exponent $p\in (2+\frac{4}{N}, 2^*)$, since, to restore compactness, it relays heavily on the uniqueness of solutions to its limited problem
$$\left\{\aligned &-\Delta u + \lambda u = |u|^{p-2}u, &x \in \mathbb{R}^{N},\\
                  &\int_{\mathbb{R}^{N}} |u|^{2}dx =c^2,          \endaligned\right.  \eqno(  S_{_{\infty},c}  )$$
and an accurate estimate on their energy levels.

\noindent\textbf{1.2 Nonlocal Case: $b>0$.} At first, we consider the limited problem(i.e., $V\equiv0$ in $(K_{V,c})$)
$$\left\{\aligned &-\left(a+ b\int_{\mathbb{R}^{N}}|\nabla u|^{2} dx  \right)\Delta u+\lambda u=|u|^{p-2}u, &x \in \mathbb{R}^{N},\\
                  &\int_{\mathbb{R}^{N}} |u|^{2}dx =c^2          \endaligned\right.  \eqno(  K_{\infty,c}  )$$
where $a, b,c>0$ are constants, $N \geq 1 $, and $ p \in (2,2^*$). The energy functional is
$$ I_{\infty}(u)=\frac{a}{2}\| \nabla u \|_2^2 + \frac{b}{4} \| \nabla u\|_2^4-\frac{1}{p}\| u \|_p^p. $$
If $N\leq 3$, by considering the minimizing problem
\begin{equation}\label{Mininonlocal}
  e(c):= \inf_{S_c} I_{V}(u),
\end{equation}
Ye \cite{Ye2014} proved that $e(c)>-\infty$ for all $c>0$ if $p\in(2,2+\frac{8}{N})$, and obtained the sharp existence of global constraint minimizer. For the case $p \in (2+\frac{8}{N},2^*)$, $e(c)=-\infty$ for all $c>0$, Ye proved the existence of one normalized solution via a suitable submanifold of $S_c$. In another word, the $L^2$-critical exponent to Kirchhoff equation is $2+\frac{8}{N}$  while it is $2+\frac{4}{N}$ to  the ``local" equation.
After that, for $2 < p < 2+\frac{4}{N}$, Zeng and Zhang \cite{Zeng2017} proved the existence and uniqueness of the minimizer to (\ref{Mininonlocal}) for any $c>0$, while for $2 +\frac{4}{N} < p < 2+\frac{8}{N}$ the authors proved that there exists a threshold
mass $c_{*} > 0$ such that for any $c \in (0, c_*)$ there is no minimizer and for $c > c_* $ there is a
unique minimizer.
Moreover, a precise formula for the minimizer and the threshold value $c_{*}$
is given according to the mass $c$.
In a recent paper \cite{Qi2022}, Qi and Zou first obtained the exact number and expressions of the positive
normalized solutions to $(K_{\infty,c})$ for $2<p\leq 2^*$, and then answered an open problem about the exact number of positive solutions of the Kirchhoff equation with
fixed frequency.
In particular, with trapping
potential $V$ and $p=2+\frac{8}{N}$, Hu and Tang \cite{Hu} considered the concentration behavior and local uniqueness
of the normalized solution for mass critical Kirchhoff equations as $a$ tending to 0. Additionally, Chen and Tang \cite{ChenAMO} considered the existence of normalized ground state solutions with $|u|^{p-2}u$ replaced by $K(x)f(u)$ where $K(x)\in\mathcal{C}(\mathbb{R}^3,\mathbb{R}^+)$ and $f(u)$ is $L^2$-supercritical. We point out that they found the minimizer on a suitable manifold following \cite{CST2021}, and the compactness was provided by the trapping nature of $K(x)$. For other related and similar results we refer readers to \cite{Li2021,Liyuhua2018,Ye2016,Zeng2021,Zhang2022}.

To the authors' knowledge, there are few result in studying bound state normalized solutions of such nonlocal problem with nonnegative potential. In the present paper, we first obtain a bound state normalized solution, which can be seen as improvements of some known results in the literature.

Indeed, new difficulties arise due to the simultaneous occurrence of nonnegative potential and nonlocal term.
As we mentioned before, the classical method in \cite{Jeanjean1997,Molle2021,Zhong2021} does not work in our situation. Thus, we adopt the linking structure constructed by \cite{Bartsch2021}. However, it seems impossible to give the accurate estimates associated with the limit problem $(K_{\infty,c})$, which
prevent  us from restoring the compactness by using the argument in \cite{Bartsch2021} directly. To overcome this obstacle we use a more subtle analysis to revisit $(K_{\infty,c})$, and find a new energy inequality and some bifurcation phenomena. With the help of a variant of the splitting lemma, we conclude our conclusions.

Our first result can be stated as follows.

\begin{Theorem}
(1)\ Assume that $1 \leq N \leq 3$, $p\in[2+\frac{8}{N},2^*)$, then $(K_{\infty,c})$ has a unique solution $(u_c,\lambda_c)$ if $c > c^*$, and has no solution if $c\leq c^*$. Moreover,
$$\| \nabla u_c \|_2^2,\ \lambda_c\rightarrow  \left\{ \aligned &0, \ \  c \rightarrow \infty,\\
 &\infty, \ \   c \rightarrow c^*, \endaligned \right.\ \ \
c^*=\left\{ \aligned &0,\qquad {\rm if} \  \   {\frac{2N+8}{N}}< p<2^*,\\
 &\left( \frac{b}{2} \right)^{\frac{N}{8-2N}} \|Q_{2+\frac{8}{N}}   \|_2^{\frac{8}{8-2N}},\  {\rm if} \  \ p= {\frac{2N+8}{N}}. \endaligned \right.$$
\noindent(2)\ Assume that $N\geq1$, $2+\frac{4}{N} <p< \min\{2+\frac{8}{N},2^*\}$, then $(K_{\infty,c})$ has a unique solution if $c=c_{1}$, no solution if $c < c_1$, and exact two solutions if $c>c_1$ where $c_1$ is a constant defined in Lemma 2.1. Moreover, if $c>c_1$, then
$$\| \nabla u_{c,1} \|_2^2\rightarrow  \left\{ \aligned &\infty, \ \  c \rightarrow \infty,\\
 & \Upsilon , \ \   c \rightarrow c_1, \endaligned \right.\ \ \
\| \nabla u_{c,2} \|_2^2\rightarrow  \left\{ \aligned &0, \ \  c \rightarrow \infty,\\
 & \Upsilon , \ \   c \rightarrow c_1, \endaligned \right.
$$
$$\lambda_{c,1}\rightarrow  \left\{ \aligned &\infty, \ \  c \rightarrow \infty,\\
 &  \Lambda, \ \   c \rightarrow c_1, \endaligned \right.\ \ \
\lambda_{c,2}\rightarrow  \left\{ \aligned &0, \ \  c \rightarrow \infty,\\
 & \Lambda , \ \   c \rightarrow c_1, \endaligned \right.
$$
where $\Upsilon$ and $\Lambda$ are constants given in Lemma 2.3.
\end{Theorem}

The second part of this paper is devoted to consider the normalized solutions to $(K_{V,c})$. Setting $m_c=I_{\infty}(u_c)$ and $u_c$ is the solution obtained in \emph{Theorem 1.1 (1)}, we introduce the following  assumptions on the potential $V$.

\noindent $(V1)$\ $1 \leq N \leq 3$, $V$ and the map $W:x \mapsto V(x)|x|$ are in $L^{\infty}(\mathbb{R}^N)$, $V \geq 0$, $\lim_{|x|\rightarrow\infty}V(x)=0$,
\begin{equation}\label{V11}
  0<\| V\|_{\infty}< \frac{2\mu}{c^2}m_c,\ \  \mu=\min\left\{1,\frac{2}{N}\right\},
\end{equation}
\begin{equation}\label{V12}
\|W\|_\infty\le \left\{  \aligned & \frac{m_c^{1/2}}{c}\left(\frac{4a(p-6)}{(p-2)^3+4(p-2)}\right)^{1/2}, \ \ N=1,\\
          & \frac{m_c^{1/2}}{c}\left(  \frac{a(N(p-2)-4)[p(2-N)+2N]^2}{4(p-2)(N(p-2)+(p(2-N)+2N))}   \right)^{1/2}, \ \ N=2,3.\endaligned\right.
\end{equation}

\noindent$(V2)$ $N=3$, $V \in L^{\frac{3}{2}}(\mathbb{R}^3)$, $W \in L^3(\mathbb{R}^3)$, $V\geq0$,
\begin{equation}\label{V21}
  \| V\|_{\frac{3}{2}}< 2aS^{2}\frac{3p-10}{9p-10},
\end{equation}
and
\begin{equation}\label{V22}
4S\| W \|_{3}\left[ \frac{3(p-2)^2}{6-p} + 1\right] + \|V \|_{\frac{3}{2}}S^2[9(p-2)+6] \leq a(3p-10)
\end{equation}
where $S$ is the best constant in the Sobolev embedding $H^1(\mathbb{R}^3) \hookrightarrow L^6(\mathbb{R}^3)$.

Notice that assumption $(V2)$ independent of $c$, contrarily to $(V1)$. Moreover, assumption $(V2)$ allows that the potential has poles, which is important for physical reasons.
We have the following results.

\begin{Theorem}
Assume $1\leq N\leq 3$, $2+\frac{8}{N} \leq p < 2^*$, and $(V1)$ holds. Then\\
(1) for every $c>c*$, there exist constant $K_1(a,N,p,c,V)$ such that
$(K_{V,c})$ has a solution $(v,\lambda) \in H^1(\mathbb{R}^N) \times \mathbb{R}^+$ if $b<K_1(a,N,p,c,V)$.\\
(2) for every $b>0$, there exist constant $K_2(a,N,p,b)$ such that
$(K_{V,c})$ has a solution $(v,\lambda) \in H^1(\mathbb{R}^N) \times \mathbb{R}^+$ if $c>K_2(a,N,p,b)$.\\
\end{Theorem}

\begin{Theorem}
Assume $N=3$ and $(V2)$ holds. Then\\
(1) for every $c>c*$, there exist constant $K_1^*(a,p,V,c)$ such that
$(K_{V,c})$ has a solution $(v,\lambda) \in H^1(\mathbb{R}^N) \times \mathbb{R}^+$ if $b<K_1^*(a,p,V,c)$.\\
(2) for every $b>0$, there exist constant $K_2^*(a,p,V,b)$ such that
$(K_{V,c})$ has a solution $(v,\lambda) \in H^1(\mathbb{R}^N) \times \mathbb{R}^+$ if $c>K_2^*(a,p,V,b)$.\\
\end{Theorem}

If $V$ satisfies the some additional regularity and monotonicity condition, then we have
\begin{Theorem}
Let $1\leq N\leq 3$, $2+\frac{8}{N} \leq p < 2^*$ and $c>c^*$. If V satisfies $(V1)$ and\\
\noindent$(V3)$
\begin{minipage}[t]{\linewidth}\ For a.e. $x \in \mathbb{R}^N$, $\nabla V(x)$ exists, $\langle \nabla V(x) \cdot x \rangle \leq 0$ and the following  map is non-decreasing
 $$
 t \mapsto t^{-5} \langle \nabla V(x/t) \cdot x \rangle, \quad t\in (0,\infty),
 $$
\end{minipage}
 then $(K_{V,c})$ has a solution $(v,\lambda) \in H^1(\mathbb{R}^N) \times \mathbb{R}^+$.
\end{Theorem}

In fact, $(V3)$ is devoted to recover the lack of compactness without the estimate in \emph{Lemma 2.4}, so the limitations on $b$ and $c$ can be dropped.
Similarly, if the potential function $V$ is radial, then our existence results can be extended to $2+\frac{4}{N}<p<2^*$ by a mountain pass argument proposed by \cite{Jeanjean1997}.

Setting $m_{c,1}=I_{\infty}(u_{c,1})$ and $m_{c,2}=I_{\infty}(u_{c,2})$, we introduce the radial assumptions on $V$

\noindent $(V4)$
\begin{minipage}[t]{\linewidth}
  $N \geq 2$, $V$ and the map $W:x \mapsto V(x)|x|$ are in $L_{rad}^{\infty}(\mathbb{R}^N)$, $V \geq 0$, $\lim_{|x|\rightarrow\infty}V(x)=0$, and \\(\ref{V11}), (\ref{V12}) hold.
\end{minipage}

\noindent $(V4')$
\begin{minipage}[t]{\linewidth}
  $N =3 $, $0\leq V \in L_{rad}^{\frac{3}{2}}(\mathbb{R}^3)$, $W\in L_{rad}^3(\mathbb{R}^3)$, and $(\ref{V22})$ holds.
\end{minipage}

\noindent $(V5)$
\begin{minipage}[t]{\linewidth}
  $N \geq 2 $, $V(x)\in L_{rad}^{\infty}(\mathbb{R}^N)$, $\nabla V(x)$ exists for a.e. $x \in \mathbb{R}^N$, $V\geq 0$, $\lim_{|x|\rightarrow\infty}V(x)=0$,
   $$\|V \|_{\infty} \leq 2c^{-2}(m_{c,2}-m_{c,1})\ \ \text{\ and\ }\  \ \langle \nabla V(x) \cdot x \rangle \leq  -\frac{N(p-2)}{p}V(x).$$
\end{minipage}

\begin{Theorem}
Let $N \geq 2$, $2+\frac{8}{N} \leq p < 2^*$, $c>c^*$, and V satisfies $(V4)$ or $(V4')$, then $(K_{V,c})$ has a solution $(v,\lambda) \in H_{rad}^1(\mathbb{R}^N) \times \mathbb{R}^+$.
\end{Theorem}
\begin{Theorem}
Let $N \geq 2$, $2+\frac{4}{N} < p < \min\{2+\frac{8}{N},2^*\}$, $c>c_1$, and V satisfies $(V5)$, then $(K_{V,c})$ has a solution $(v,\lambda) \in H_{rad}^1(\mathbb{R}^N) \times \mathbb{R}^+$.
\end{Theorem}

This paper is structured as follows. In Section 2, we address the study of the limit problem and prove \emph{Theorem 1.1}. Afterwards, in Section 3, we present the linking geometry which is crucial for the proof of our results. In Section 4, we prove \emph{Theorem 1.2}, \emph{Theorem 1.3} and \emph{Theorem 1.4}. Finally, in Section 5 we consider the case of radial potentials and conclude the proof of \emph{Theorem 1.5} and \emph{Theorem 1.6}.

\noindent \textbf{Notations:} From now on in this paper, otherwise mentioned, we use the following notations:

$\bullet$ $L^p(\mathbb{R}^N)$ with $ p \in [1,\infty)$ is the Lebesgue space with the norm $\|u \|_p=(\int_{\mathbb{R}^N}|u|^p dx )^{1/p}$.

$\bullet$ $L^{\infty}(\mathbb{R}^N)$ is the Lebesgue space with the norm $\|u \|_{\infty}=\inf\{C: |u(x)|\leq C \text{\ a.e.\ in\ } \mathbb{R}^N\}$.

$\bullet$ $H^{1}(\mathbb{R}^{N})$ is the usual Sobolev space with the norm
$\| u \|_{H^{1}(\mathbb{R}^N)}=\left( \int_{\mathbb{R}^N}|\nabla u|^2+|u|^2 dx\right)^{1/2}$.

$\bullet$ $H^1_{rad}(\mathbb{R}^{N})=\{ u \in H^{1}(\mathbb{R}^{N}):u\ \text{is\ radial\ symmetric}\}$.

$\bullet$ $L_{rad}^p(\mathbb{R}^N)$ with $ p \in [1,\infty]$ denotes $\{ u \in L^p(\mathbb{R}^N):u\ \text{is\ radial\ symmetric}\}$.

$\bullet$ $t \star u(\cdot) = t^{\frac{N}{2}}u(t\cdot)$ for all $u \in H^{1}(\mathbb{R}^N)$ and $t>0$.

$\bullet$
$a\sim b$ means that there exist constants $0<C'\leq C$ such that $C'a \leq b \leq Ca$.

\section{The limited problem}

In this section, we study the limited problem

$$\left\{\aligned &-\left(a+ b\int_{\mathbb{R}^{N}}|\nabla u|^{2} dx  \right)\Delta u+\lambda u=|u|^{p-2}u, &x \in \mathbb{R}^{N},\\
                  &\int_{\mathbb{R}^{N}} |u|^{2}dx =c^2          \endaligned\right.  \eqno(  K_{\infty,c}  )$$
where $a, b,c>0$ are constants, $N \geq 1 $, and $ p \in (2,2^*$).
The normalized solutions to $(  K_{\infty,c}  )$ are obtained by looking  for critical points of $\mathcal{C}^1$ functional
 \begin{equation}\label{equ:E F infty}
   I_{\infty}(u)=\frac{a}{2}\| \nabla u \|_2^2 + \frac{b}{4} \| \nabla u\|_2^4-\frac{1}{p}\| u \|_p^p
 \end{equation}
constrained on the $ L^2$-sphere in $ H^1(\mathbb{R}^{N})$:
$$   S_c = \{  u \in H^1(\mathbb{R}^N) \ | \    \|u \|_2^2 = c^2  >0    \}.         $$

If $(u_c, \lambda_c)$ is a solution to $(K_{\infty,c})$ for suitable $c>0$, then $u_c$ satisfies the following Pohozaev and Nehari identity
\begin{equation}\label{equ: pohozaev infty}
 P_{\infty}(u_c):= a\| \nabla u_c \|_2^2 + b\| \nabla u_c \|_2^4- \frac{N(p-2)}{2p} \|u_c \|_p^p=0,
\end{equation}
\begin{equation}\label{equ:nehari infty}
 N_{\infty,\lambda_c}(u_c):= a\| \nabla u_c \|_2^2 + b\| \nabla u_c \|_2^4 + \lambda_c \|u_c\|_2^2-\|u_c \|_p^p=0.
\end{equation}
Let $D:=\left(\int_{\mathbb{R}^N}|\nabla u_c|^2dx\right)^{1/2}$, then
\begin{equation}\label{equ: lambda c}
  \lambda_c =\frac{2N-p(N-2)}{N(p-2)c^2}D^2\left( a+bD^2 \right).
\end{equation}
Therefore, we see that $u_c \in S_c$ is a positive solution of the following equation
$$ -\frac{N(p-2)}{4} \triangle u_c + \frac{2N-p(N-2)}{4} \frac{D^2}{c^2}u_c = \frac{N(p-2)}{4(a+bD^2)} |u_c|^{p-2}u_c,\ \ x \in \mathbb{R}^{N}. $$
In spirit of \cite[Lemma 2.2]{YeLi2019}, by the uniqueness of positive solutions(up to translations) to (\ref{equ:Q sovles equ}) and a suitable rescaling, we have
\begin{equation}\label{precise solution infty}
u_c(x) = \left[ \frac{4a+4bD^2}{N(p-2)}   \right]^{\frac{1}{p-2}} \left(  \frac{D}{c}  \right)^{\frac{2}{p-2}} Q_p \left(\frac{D}{c}x  \right),
\end{equation}
and
$$   \left[ \frac{4a+4bD^2}{N(p-2)}   \right]^{\frac{2}{p-2}} \left(  \frac{D}{c}  \right)^{\frac{4-N(p-2)}{p-2}} \| Q_p   \|_2^2 =  c^2  = \| u_c\|_2^2    .       $$
More precisely, $D$ satisfies
\begin{equation}\label{equ:Root equation}
  D^2\left[  \frac{4}{N(p-2)}  \right]^\frac{4}{4-N(p-2)}\| Q_p \|_2^{\frac{2(p-2)}{4-N(p-2)}} = c^{\frac{4N-2p(N-2)}{4-N(p-2)}}\left( a+bD^2 \right)^\frac{4}{N(p-2)-4}.
\end{equation}
In this way, $(K_{\infty, c} )$ has $\boldsymbol{k}$ distinct solutions if and only if (\ref{equ:Root equation})(treat $D>0$ as an unknown in the equation) has $\boldsymbol{k}$ distinct roots. Hence, we can obtain the exact number of solutions and the precise characterization by studying (\ref{equ:Root equation}).

Then the first result is as follows.

\begin{Lem} \label{Lem:existence}
Assume that $N \geq 1$, $2+\frac{4}{N}<p<2^*$, and $\theta=N(p-2)-4 $, $ \eta=8-N(p-2)$.\\
(1) \begin{minipage}[t]{\linewidth}
     If $N \leq 3$, $p \in (2+\frac{8}{N},2^*)$, then $(K_{\infty,c})$ has a unique solution for $c>0$.
    \end{minipage}
(2)\ \begin{minipage}[t]{\linewidth}
     If $N \leq 3$, $p=2+\frac{8}{N}$, then $(K_{\infty,c})$ has a unique solution for $c>c_0$, no solution for $c\leq c_0$.
    \end{minipage}\\
(3)
    \begin{minipage}[t]{\linewidth}
     If $1 \leq N \leq 4$, $p\in (2+\frac{4}{N},2^*)$, then $(K_{\infty,c})$ has a unique solution for $c=c_{1}$, no solution\\ for $c < c_1$, and exact two solutions for $c>c_1$.
    \end{minipage} \\
Moreover, $$c_0=  \left( \frac{b}{2} \right)^{\frac{N}{8-2N}} \|Q_{2+\frac{8}{N}}   \|_2^{\frac{8}{8-2N}} ,$$ $$   c_1=  \left(\frac{16}{N(p-2)}\right)^{\frac{2}{2N-p(N-2)}} \|Q_p \|_2^{\frac{2(p-2)}{2N-p(N-2)}}      \left(    \frac{b}{\theta}  \right)^{\frac{\theta}{4N-2p(N-2)}}      \left(  \frac{a}{\eta}  \right)^{\frac{\eta}{4N-2p(N-2)}} ,$$
and the solutions can be characterized as
$$ u_{c,i}=\left[ \frac{4a+4bD_i^2}{N(p-2)}   \right]^{\frac{1}{p-2}} \left(  \frac{D_i}{c}  \right)^{\frac{2}{p-2}} Q_p \left(\frac{D_i}{c}x  \right), $$
$$ \lambda_{c,i} =\frac{2N-p(N-2)}{N(p-2)c^2}D_i^2\left( a+bD_i^2 \right), i=1,2, $$
where $D_i$ are the roots of (\ref{equ:Root equation}).

In particular, $(K_{\infty,c})$ has unique solution $(u_c,\lambda_c)$ if $D_1 = D_2$.

\end{Lem}
\noindent \textbf{Proof:}\
From direct calculation, (\ref{equ: lambda c}) and (\ref{precise solution infty}), the conclusion is obvious.  \qed

Next, we give the ``variational characterization" to those solutions obtained above.

\begin{Lem}

Assume that $p \in [2+\frac{8}{N},2^*)$, $u_c$ is the unique solution to $(K_{\infty,c})$. Then,

\noindent(1)\ \ \begin{minipage}[t]{\linewidth}
     $u_c$ is mountain-pass type solution. That is, there exists a constant $r(c)>0$ such that
$$ I_{\infty}(u_c)=m_c=\inf_{\sigma \in \Sigma_c} \max_{t\in[0,1]} I_{\infty}(\sigma(t))\ \ \ \ \  $$
where
$$  \Sigma_c =\left\{ \sigma \in \mathcal{C}([0,1],S_c)~:~\sigma(0) \in A_{r(c)}   ,~I_{\infty}(\sigma(1))   < 0                         \right\} $$
   and
   $$ A_{r(c)}=\{u\in {S}_c : \| \nabla u\|_2^2 \leq r(c)\}.$$
\end{minipage}

\noindent(2) \begin{minipage}[t]{\linewidth}
      The maps $c \mapsto m_c, \| \nabla u_c \|_2^2 $, $\lambda_c$ are continuous and strictly decreasing on $(c^*, \infty)$.  \end{minipage}
\end{Lem}
\noindent(3) \begin{minipage}[t]{\linewidth}

$$\| \nabla u_c \|_2^2,\ m_c,\ \lambda_c,\ \frac{m_c}{c^2} \rightarrow  \left\{ \aligned &0, \ \  c \rightarrow \infty,\ \ \ \ {\frac{2N+8}{N}}\leq p<2^*,\\
&\infty,\ \   c \rightarrow 0,\ \ \ \     {\frac{2N+8}{N}}< p<2^*, \\
 &\infty, \ \   c \rightarrow c_0, \ \ \ \ p= {\frac{2N+8}{N}}. \endaligned \right.$$
\end{minipage}
\textbf{Proof:}\ (1)\ Similar to\cite{Jeanjean1997} and \cite{Ye2015}, we omit it.\\
\noindent(2)\
From (\ref{equ: pohozaev infty}) and (\ref{equ: lambda c}), we have
$$  \aligned m_c=I_{\infty}(u_c)&= \frac{N(p-2)-4}{2N(p-2)}a\| \nabla u_c \|_2^2 + \frac{N(p-2)-8}{4N(p-2)}b\| \nabla u_c  \|_2^4    \endaligned $$
and
$$ \lambda_c =\frac{\zeta}{N(p-2)}  \left(\frac{N(p-2)}{4}\right)^{\frac{4}{\zeta}}\| Q_p \|_2^{\frac{2-p}{\zeta}}(\| \nabla u_c \|_2^2)^{\frac{2(p-2)}{\zeta}} \left( a+b\| \nabla u_c \|_2^2 \right)^\frac{(p-2)(2-N)}{\zeta}.       $$
where $\zeta=2N-p(N-2)$. By the property of (\ref{equ:Root equation}), $c\mapsto \|\nabla u_c \|_2^2$ is strictly decreasing.
So the proof is complete.\\
\noindent(3)\ From (\ref{equ:Root equation}), if $ 2+\frac{8}{N} <p<2^*$, we have
$$ \|\nabla u_c \|_2^2 \sim (c^2)^\frac{2N-p(N-2)}{8-N(p-2)}.    $$
Then,
$$ \frac{m_c}{c^2} \sim \frac{N(p-2)-4}{2N(p-2)} a (c^2)^{\frac{2p-8}{8-N(p-2)}}+ \frac{N(p-2)-8}{4N(p-2)} b (c^2)^{\frac{(4-N)(p-2)}{8-N(p-2)}}. $$
In particular, if $p=2+\frac{8}{N}$, we have
$$  \|\nabla u_c \|_2^2 = \frac{c^{\frac{2N-8}{N}}a\| Q_{\frac{2N+8}{N}} \|_2^{\frac{4}{N}}  }{2-b\| Q_{\frac{2N+8}{N}} \|_2^{\frac{4}{N}} c^{\frac{2N-8}{N}}    }.     $$
Then,
$$  \frac{m_c}{c^2} = \frac{a}{4}  \frac{c^{-\frac{8}{N}}a\| Q_{\frac{2N+8}{N}} \|_2^{\frac{4}{N}}  }{2-b\| Q_{\frac{2N+8}{N}} \|_2^{\frac{4}{N}} c^{\frac{2N-8}{N}}    }.   $$
Combining (\ref{equ: lambda c}) and direct calculation, the proof is complete.  \qed

\begin{Lem}
Assume that $N\geq2$, $p\in\left(\frac{2N+4}{N}, \min\{\frac{2N+8}{N},2^*\}\right)$, $c>c_1$, and $D_1 > D_2$. Then,\\
(1)\begin{minipage}[t]{\linewidth}
\ \  $u_{c,1}$ is a local minimizer, and $u_{c,2}$ is a mountain-pass solution. That is,
$$ m_{c,1} =I_{\infty}(u_{c,1}) = \inf_{v\in {S}_{c}^r \setminus A_{T}^r} I_{\infty}(v),                      $$
$$ m_{c,2} =I_{\infty}(u_{c,2}) = \inf_{\gamma \in \Gamma_c} \max_{s \in [0,1]} I_{\infty}(\gamma (s) )                                        $$
where
$$ {S}_{c}^r=\{ u\in H_{rad}^{1}(\mathbb{R}^N):\|u \|_2^2 = c^2  \}, \quad  A_{T}^r=\{ u \in S_{c}^r:\| \nabla u\|_2^2 \leq T \}  $$
$$  \Gamma_c = \{ \mathcal{C}([0,1],{S}_{c}^r): \gamma(0)\in A_{l}^r,\gamma(1)=u_{c,1} \}, \quad  A_{l}^r=\{ u \in S_{c}^r:\| \nabla u\|_2^2 \leq l \}, $$
 and $T,l$ are constants depending on c.
\end{minipage}
\noindent(2) \begin{minipage}[t]{\linewidth}There exist constants $\Upsilon$ and $\Lambda$ depending on c such that
      $$\| \nabla u_{c,1} \|_2^2\rightarrow  \left\{ \aligned &\infty, \ \  c \rightarrow \infty,\\
 & \Upsilon^+ , \ \   c \rightarrow c_1, \endaligned \right.\ \ \
\| \nabla u_{c,2} \|_2^2\rightarrow  \left\{ \aligned &0, \ \  c \rightarrow \infty,\\
 & \Upsilon^- , \ \   c \rightarrow c_1, \endaligned \right.
$$
$$\lambda_{c,1}\rightarrow  \left\{ \aligned &\infty, \ \  c \rightarrow \infty,\\
 &  \Lambda, \ \   c \rightarrow c_1, \endaligned \right.\ \ \
\lambda_{c,2}\rightarrow  \left\{ \aligned &0, \ \  c \rightarrow \infty,\\
 & \Lambda , \ \   c \rightarrow c_1. \endaligned \right.
$$
Moreover,
$$\left\{ \aligned & m_{c,1} \longrightarrow  m_{c_1}^-\\
             & m_{c,2} \longrightarrow m_{c_1}^+ \endaligned\right. \quad \text{as}\ c\rightarrow c_1 \quad  \ \text{and}\ \quad  \left\{\aligned & m_{c,1} \longrightarrow  -\infty\\
                                                                                                            & m_{c,2} \longrightarrow 0    \endaligned\right. \quad \text{as}\ c\rightarrow \infty. $$

    \end{minipage}
\end{Lem}
\noindent\textbf{Proof:}\ (1) From (\ref{c GN ineq}), we have
\begin{equation}\label{Bar I}
 \aligned I_{\infty}(u)& = \frac{a}{2}\| \nabla u \|_{2}^2 + \frac{b}{4}\|\nabla u \|_2^4 -\frac{1}{p} \| u \|_{p}^{p}  \\
                      &\geq  \frac{a}{2}\| \nabla u \|_{2}^2 + \frac{b}{4}\|\nabla u \|_2^4 - \frac{c^{p-\frac{N(p-2)}{2}}}{2\| Q_p \|_2^{p-2}}  \| \nabla u \|_2^{\frac{N(p-2)}{2}}\\
                      &:= \bar{I}_{\infty}(u),\ \ \forall u \in S_{c}^r, \endaligned
\end{equation}
with equality only for $u=t \star \left(\frac{c}{\|Q_p\|_2}Q_p\right)$ for all $t\in \mathbb{R}$, up to translations.
\noindent Set auxiliary function $f(t):[0, \infty) \mapsto \mathbb{R}$ by $f(t)=  \frac{a}{2}t + \frac{b}{4}t^2 - \frac{c^{p-\frac{N(p-2)}{2}}}{2\| Q_p \|_2^{p-2}}  t^{\frac{N(p-2)}{4}}$, and it is easy to see that the nonzero critical points of $f$ satisfy
$$   t \left[  \frac{4}{N(p-2)}  \right]^\frac{4}{4-N(p-2)}\| Q_p \|_2^{\frac{2(p-2)}{4-N(p-2)}} = c^{\frac{4N-2p(N-2)}{4-N(p-2)}}\left( a+bt \right)^\frac{4}{N(p-2)-4}.       $$
We may assume $t_1 > t_2>0$ such that $f'(t_1)=f'(t_2)=0$, $f''(t_1)>0$, $f''(t_2)<0$, then there exist constants $T$, $r$ depending on c such that $f'(T)<0$, $  t_2 < T < t_1-r  $, $f(t_1)=\inf\limits_{t \geq T}f(t)$, and $\inf\limits_{|t-t_1|=r} f(t) > f(t_1)$.

Correspondingly, we have \\
\begin{equation}\label{equ: end of mp}
  \begin{split}
&D_2^2 < T < D_1^2  ,\ \ \ I_{\infty}(u_{c,1})=\bar{I}_{\infty}(u_{c,1})=\inf\limits_{S_{c}^r \setminus A_T^r}\bar{I}_{\infty}(u_{c,1}),\\
&\inf\limits_{ \left|\|\nabla u\|_2^2 -\|\nabla u_{c,1}\|_2^2\right| =r         }\bar{I}_{\infty}(u) >\bar{I}_{\infty}(u_{c,1}).
\end{split}
\end{equation}
Noticed that $ \inf\limits_{{S}_{c}^r \setminus A_T^r} I_{\infty}(v) \geq \inf\limits_{S_{c}^r \setminus A_T^r}\bar{I}_{\infty}(u_{c,1})$ from (\ref{Bar I}), and $t_1$ is local minimizer of $f(t)$, we get the conclusion.

Similar to \cite{Jeanjean1997}, there exist $0 < l< 2l <T$ such that
$$  0< \sup\limits_{u\in A_{l}^r} I_{\infty}(u)<\inf\limits_{u\in\partial A_{2l}^r}I_{\infty}(u)\ \ \text{and}\ \ I_{\infty}(u)>0\ \ \text{for all}\ u\in A_{l}^r.  $$
Then for any $\gamma\in\Gamma_c$($\neq\emptyset$ which is obvious),there exist $s_0,s_1>0$ such that \
$$  \| \gamma(s_0) \|_2^2=l \ \ \text{and}\ \ \left| \|\gamma(s_1) \|_2^2 -\|  u_{c,1} \|_2^2      \right| =r,      $$
So combining this with (\ref{equ: end of mp}), we obtain $\inf_{\gamma \in \Gamma_c} \max_{s \in [0,1]} I_{\infty}(\gamma (s) ) > \max\{ \gamma(0),\gamma(1)   \}.$

Following \cite{Jeanjean1997}, there exist $\{ u_n\} \subset S_{c}^r$ such that
$$P_{\infty}(u_n) \rightarrow 0 \text{\ \ and\ \ } \|  I'_{\infty}|_{S_{c}^r} (u_n) \| \rightarrow 0.$$
Since $H_{r}^1(\mathbb{R}^N) \hookrightarrow L^p(\mathbb{R}^N) $ is compact for $p\in(2,2^*)$, then there exists $(\tilde{u},\tilde{\lambda}) \in H^{1}(\mathbb{R}^N)\setminus{0} \times \mathbb{R}^+$
solves $(K_{\infty,c})$ by standard argument. From \emph{Lemma 2.1}, we deduce $(\tilde{u},\tilde{\lambda})=(u_{c,2},\lambda_{c,2})$.

\noindent(2) Set $\Upsilon$ is the unique root of (\ref{equ:Root equation}) in the case of $c=c_1$ and $\Lambda=\frac{2N-p(N-2)}{N(p-2)c_1^2}\Upsilon^2\left( a+b\Upsilon^2 \right)$, then
from the property of (\ref{equ:Root equation})(or $f(t)$) which implies $D_1^2 \rightarrow 0$, $D_2^2 \rightarrow \infty$ as $c\rightarrow \infty$ and $D_1^2, D_2^2 \rightarrow \Upsilon$ as $c \rightarrow c_1$, we prove the first part.

From (\ref{equ: pohozaev infty}) and \emph{Lemma 2.1}, we get
$$ m_{c,i}=I_{\infty}(u_{c,i}) = \frac{N(p-2)-4}{2N(p-2)} a D_i^2+ \frac{N(p-2)-8}{4N(p-2)} b D_i^4.  $$
Since $\frac{N(p-2)-8}{4N(p-2)} < 0$, we prove the second part.  \qed

\noindent\textbf{Proof of Theorem 1.1}:\ Combining \emph{Lemma 2.1}, \emph{Lemma 2.2} with \emph{Lemma 2.3}, we deduce the conclusion. \qed


In order to show the connection between nonlocal case $(b>0)$ and local case $(b=0)$, we give a estimate of the energy level and study the asymptotic profiles of solutions when $b \rightarrow 0$. We would like to emphasize that the following lemmas are helpful to recover the compactness.

\begin{Lem}
If $2+\frac{8}{N} \leq p < 2^*$, $\alpha \geq \beta > c_0$, then $  \frac{m_\beta}{m_\alpha} \geq \frac{\alpha^q}{\beta^q} $ with $q=\frac{4N-2p(N-2)}{N(p-2)-4}$.
\end{Lem}
\noindent \textbf{Proof:}\ Defining a auxiliary function $$f(t)= \frac{N(p-2)-4}{2N(p-2)}a(a+bt)^{\frac{4}{N(p-2)-4}}+\frac{N(p-2)-8}{4N(p-2)}bt(a+bt)^{\frac{4}{N(p-2)-4}}$$
for $2+\frac{8}{N} \leq p <2^*$. It is standard to see $f(t)$ is non-decreasing. Set $D_{\alpha}^2=\| \nabla u_{\alpha} \|_2^2 $ and $D_{\beta}^2 = \| \nabla u_\beta \|_2^2$, then $D_{\alpha}^2 \leq D_{\beta}^2$ due to the above Lemma. So, we get $ f(D_{\alpha}^2) \leq f(D_{\beta}^2)$ i.e.
$$\aligned
         & \frac{N(p-2)-4}{2N(p-2)}a(a+bD_{\alpha}^2)^{\frac{4}{N(p-2)-4}}+\frac{N(p-2)-8}{4N(p-2)}bD_{\alpha}^2(a+bD_{\alpha}^2)^{\frac{4}{N(p-2)-4}} \\
    \leq & \frac{N(p-2)-4}{2N(p-2)}a(a+bD_{\beta}^2)^{\frac{4}{N(p-2)-4}}+\frac{N(p-2)-8}{4N(p-2)}bD_{\beta}^2(a+bD_{\beta}^2)^{\frac{4}{N(p-2)-4}},
\endaligned $$
which means
$$  \frac{D_{\alpha}^{-2}(a+bD_{\alpha}^2)^\frac{4}{N(p-2)-4}}{D_{\beta}^{-2}(a+bD_{\beta}^2)^\frac{4}{N(p-2)-4}}  \leq \frac{\frac{N(p-2)-4}{2N(p-2)} a D_{\beta}^2+ \frac{N(p-2)-8}{4N(p-2)} b D_{\beta}^4}{\frac{N(p-2)-4}{2N(p-2)} a D_{\alpha}^2+ \frac{N(p-2)-8}{4N(p-2)} b D_{\alpha}^4}    =\frac{I_{\infty}(u_{\beta})}{I_{\infty}(u_{\alpha})}  =\frac{m_{\beta}}{m_{\alpha}}.             $$
Combining (\ref{equ:Root equation}), we deduce $\frac{m_{\beta}}{m_{\alpha}} \geq \frac{\alpha^q}{\beta^q}.$  \qed

\begin{Lem}
Let $b \rightarrow 0$, then $(u_{c,b}, \lambda_{c,b}) \rightarrow (Z_c, \Lambda_c) $ in $H^1(\mathbb{R}^N) \times \mathbb{R}$ strongly where $Z_c$ is the unique normalized solution to
$$ -a \triangle u + \lambda u = |u|^{p-2}u,\ \ \| u \|_2 = c. $$
\end{Lem}
\noindent\textbf{Proof:}\ From (\ref{precise solution infty}) and (\ref{equ:Root equation}), we deduce the conclusion. \qed

\section{Linking Geometry}
In this section, we adopt the min-max argument developed by \cite{Bartsch2021} to $I$ and obtain linking geometry under different assumptions on $V$.
\subsection{Under $(V1)$ or $(V2)$ assumption}
To begin with, we always assume $p \in [\frac{2N+8}{N},2^*)$, $c>c^*$ which implies $(K_{\infty,c})$ has unique radial positive solution $u_c$, and $(V1)$ or $(V2)$ holds.

For $R>0$ that will be determined later, we set
$ Q=B_R \times [h_1,h_2] \subset \mathbb{R}^N \times \mathbb{R} $ where $B_R=\{x \in \mathbb{R}^N:|x|\leq R\}$ and $0 < h_1 <h_2$. If there exist a suitable choice of $Q$ such that
$$   \sup_{\phi \in \Phi_c} \max_{(y,h)\in \partial Q} I(\phi(y,h)) < m_{V,c}:= \inf_{\phi\in\Phi_c}\max_{(y,h)\in Q}I(\phi(y,h))   $$
where
$$ \Phi_c =\{  \phi\in \mathcal{C}(Q, S_c) : \phi(y,h)=h\star u_c (\cdot - y) \ \text{for all}\ (y,h)\in \partial Q            \},          $$
we can adopt the argument of \cite{Bartsch2021}. We note that the first key step is to distinguish between $m_{V,c}$ and $m_c$, and this relays heavily on the uniqueness of solution to $(K_{\infty,c})$.
Before doing this, we need to recall some notations and preliminary results.

In order to study the behavior of $PS$ sequence, we introduce a splitting lemma which plays a crucial role in obtaining the new linking geometry and overcoming the lack of compactness.
For $\lambda>0$ we set
$$I_{\lambda}(u)= \frac{a}{2}\int_{\mathbb{R}^N}|\nabla u|^2 dx + \frac{1}{2} \int_{\mathbb{R}^N}(V(x)+\lambda) u^2 dx + \frac{b}{4} \left(\int_{\mathbb{R}N} |\nabla u|^2 dx\right)^2 - \frac{1}{p}\int_{\mathbb{R}^N} |u|^p dx  $$
and
$$I_{\infty,\lambda}(u)=\frac{a}{2}\int_{\mathbb{R}^N}|\nabla u|^2 dx + \frac{1}{2} \int_{\mathbb{R}^N} \lambda u^2 dx + \frac{b}{4} \left(\int_{\mathbb{R}N} |\nabla u|^2 dx\right)^2 - \frac{1}{p}\int_{\mathbb{R}^N} |u|^p dx.$$

\begin{Prop}
Let $\{v_n\} \subset H^{1}(\mathbb{R}^N) $ be a (PS) sequence for $I_{\lambda}$ such that $v_n \rightharpoonup v$ in $H^{1}(\mathbb{R}^N)$ and $ \lim\limits_{n \rightarrow \infty}\| \nabla v_n \|_2^2 = A^2$. Then there exist an integer $k \geq 0$, $k$ non-trivial solutions $ w^1, w^2, \dots, w^k $ to the equation
\begin{equation}\label{equ:A equation}
  -(a+bA^2)\triangle w + \lambda w = |w|^{p-2}w
\end{equation}
and $k$ sequences $\left\{y_{n}^j\right\}\in H^{1}(\mathbb{R}^N)$, $1 \leq j \leq k$, such that $|y_n^j|\rightarrow \infty$ as $n \rightarrow \infty$.

Moreover, we have
\begin{equation}\label{dec equ}
\begin{split}
    & v_n - \sum_{j=1}^{k}w^{j}(\cdot - y_{n}^{j} )\rightarrow v \ \ \text{in}\ \ H^{1}(\mathbb{R}^N), \\
    & \|  v_n \|_2^2 \rightarrow \| v \|_2^2 + \sum_{j=1}^{k}\| w^j \|_2^2,\\
    & A^2=\|\nabla v \|_2^2 + \sum_{j=1}^{k} \| \nabla w^j \|_2^2,
\end{split}
\end{equation}
and
\begin{equation}\label{mp level des}
  I_{\lambda}(v_n) \rightarrow J_{\lambda}(v) + \sum_{j=1}^{k}J_{\infty}(w^j)
\end{equation}
as $ n \rightarrow \infty$ where
$$\aligned
   & J_{\lambda}(u)= \left(\frac{a}{2}+\frac{bA^2}{4}\right)\int_{\mathbb{R}^N}|\nabla u|^2 dx + \frac{1}{2} \int_{\mathbb{R}^N}(V(x)+\lambda) u^2 dx - \frac{1}{p}\int_{\mathbb{R}^N} |u|^p dx; \\
   & J_{\infty}(u)= \left(\frac{a}{2}+\frac{bA^2}{4}\right)\int_{\mathbb{R}^N}|\nabla u|^2 dx + \frac{1}{2} \int_{\mathbb{R}^N}\lambda u^2 dx - \frac{1}{p}\int_{\mathbb{R}^N} |u|^p dx.
\endaligned $$
\end{Prop}
\noindent\textbf{Proof:}\ \ See \cite[Proposition 2.1]{Xe2016}. \qed

Then recalling the notation of $barycentre$ of $u \in H^{1}(\mathbb{R}^N) \setminus \{0\}$.
Setting $$    \alpha(u)(x)=\frac{1}{|B_1(0)|} \int_{_{B_{1}(x)}} |u(y)|dy,   $$
we observe $\alpha (u)$ is bounded and continuous, so the function
$$ \hat{u}(x) =\left[     \alpha(u)(x)-\frac{1}{2}\max\alpha(u)             \right]^{+}                               $$
is well-defined, continuous, and has compact support. Therefore we can define $\beta:H^{1}(\mathbb{R}^N) \setminus \{0\} \rightarrow \mathbb{R}^N$ as
$$  \beta(u) = \frac{1}{\|  \hat{u} \|_1}\int_{\mathbb{R}^N}\hat{u}(x)xdx.   $$
The map $\beta$ is well-defined and enjoys following properties:\\
$\bullet$ \ $\beta$ is continuous in $H^1(\mathbb{R}^N) \setminus \{0\}$;  \\
$\bullet$ \ $\beta(u)=0$ for all $u \in H_{rad}^1$;  \\
$\bullet$ \ $\beta(tu)=\beta(u)$ for all $t\neq 0$ and for all $u \in H^1({\mathbb{R}^N}) \setminus \{0\}$;\\
$\bullet$ \ setting $u_z(x)=u(x-z)$ for all $z\in\mathbb{R}^N$ and $u \in H^{1}(\mathbb{R}^N)\setminus \{0\}$ there holds $\beta(u_z)=\beta(u)+z$.  \\
Now we define
$$  \aligned  &\mathcal{D}   := \{ D \subset S_c : \text{D is compact,connect},\ h_1\star U_c, h_2\star U_c \in D         \},\\
              &\mathcal{D}_0 := \{ D\in \mathcal{D}:\beta(u)=0\ \text{for all}\ u \in D       \}, \\
              &\mathcal{D}_r :=  \mathcal{D} \cap H_{rad}^{1}(\mathbb{R}^N), \endaligned   $$
and
$$ \aligned  & l_c := \inf_{D\in\mathcal{D}}\max_{u\in D}I_{\infty}(u)  \\
             & l_c^0 := \inf_{D\in\mathcal{D}_0}\max_{u\in D}I_{\infty}(u)   \\
             & l_c^r := \inf_{D\in\mathcal{D}_r}\max_{u\in D}I_{\infty}(u).                                  \endaligned     $$
\begin{Lem}
$l_c=l_c^0=l_c^r=m_c$.
\end{Lem}
\noindent\textbf{Proof:} Clearly $l_c^r \geq l_c^0 \geq l_c$ since $ \mathcal{D}_r \subset \mathcal{D}_0 \subset \mathcal{D}$. We only need to prove $
l_c \geq m_c$ and $m_c \geq l_c^r$.

Arguing by contradiction that $m_c > l_c$. Then $\max_{u \in D}I_{\infty}(u)<m_c$ for some $D \in \mathcal{D}$, hence $ \sup_{u \in V_{\delta}(D)}I_{\infty}(u)<m_c$ for some $\delta>0$; here $V_\delta(D)$ is the open $\delta$-neighborhood of $D$ and $V_{\delta}(D)$ is a path-connected. There exists a path $\sigma \in \Lambda_c$ with a suitable choice of $Q$ such that $\max_{t\in[0,1]}I_{\infty}(\sigma(t))$, a contradiction.

The equality $m_c \geq l_{c}^r$ follows from the fact that the set $D:=\{ h\star u_c : h\in [h_1,h_2] \} \in \mathcal{D}_r$ satisfies
$$\max_{D}I_{\infty}(u)=\max_{h\in[h_1,h_2]}I_{\infty}(h\star u_c)=m_c.   $$  \qed

The following Lemma is from \cite{Bartsch2021}.
\begin{Lem}
Let $E$ be a Hilbert manifold and let $F\in \mathcal{C}^1(E,\mathbb{R})$ be a given functional. Let $T \subset E$ be compact and consider a subset
$$ \mathcal{X} \subset \{  X \subset E : X\ \text{is compact},\ T \subset X  \} $$
which is homotopy-stable, i.e. it is invariant with respect to deformations leaving $T$ fixed. Assume that
$$\max_{u\in T}F(u)<c:=\inf_{X \in \mathcal{X}}\max_{u\in X}F(u) \in \mathbb{R}.$$
Let $\varepsilon_n \in \mathbb{R}$ be such that $\varepsilon \rightarrow 0$ and $X_n \in \mathcal{X}$ be a sequence such that
$$ 0 \leq \max_{u \in X_n}F(u)-c \leq \varepsilon_n. $$
Then there exists a sequence $v_n \in E$ such that \\
(1) $|F(v_n)-c|\leq\varepsilon_n$,\\
(2) $\| \nabla_E J(v_n) \| \leq \tilde{c}\sqrt{\varepsilon_n}$,\\
(3) $dist(v_n, X_n)\leq \tilde{c}\sqrt{\varepsilon_n}$, \\
for some constant $\tilde{c}$.
\end{Lem}

\begin{Lem}
$L_c:=\inf_{D \in \mathcal{D}_0}\max_{u\in D}I(u) > m_c.  $
\end{Lem}
\noindent\textbf{Proof:} Since $V \geq 0$ and $l_c^0 = m_c$, we have
\begin{equation}\label{geq side}
\max_{u\in D}I(u) \geq \max_{u \in D}I_{\infty}(u) \geq l_c^0 = m_c,  \ \ \ \text{for all}\ D \in \mathcal{D}_0.
\end{equation}
Now we argue by contradiction and assume that there exists a sequence $D_n \in \mathcal{D}_0$ such that
$$\max_{u \in D_n}I(u)\rightarrow m_c\ \ \text{and}\ \ \max_{u\in D_n}I_{\infty}(u)\rightarrow m_c. $$
We define the functional $\tilde{I}_{\infty}:H^1(\mathbb{R}^N)\times \mathbb{R} \rightarrow \mathbb{R} $ by $\tilde{I}_{\infty}(u,h):=I_{\infty}(h\star u)$ constrained to $E:=S_c\times\mathbb{R}$, and apply Lemma 3.3 with
$$ T:=\{ (h_1\star u_c, 1), (h_2\star u_c, 1)   \}       $$
and
$$ \mathcal{X}:=\{ X\in E : X\ \text{is compact, connected,}\ T\subset X   \}.   $$
Note that
$$ \tilde{l}_c := \inf_{X \in \mathcal{X}}\max_{(u,h)\in C}\tilde{I}_{\infty}(u,h)=l_c=m_c    $$
since $\mathcal{D}\times\{1\} \subset \mathcal{X}$, hence $l_c \geq \tilde{l}_c$, and for any $X\in \mathcal{X}$ we have $D:=\{ h\star u:(u,h) \in X \}\in \mathcal{D}$
and
$$  \max_{(u,h)\in X}\tilde{I}_{\infty}(u,h)=\max_{(u,h)\in X}I_{\infty}(h\star u)=\max_{v \in D}I_{\infty}(v)  $$
hence $ l_c \leq \tilde{l}_c $. Hence, Lemma 3.3 yields a sequence $(u_n,h_n) \in S_c \times \mathbb{R}$ such that\\
(1) $|\tilde{I}_{\infty}(u_n,h_n)-m_c|\rightarrow 0\ \text{as}\ n\rightarrow \infty,$\\
(2) $\| \nabla_{S_c \times \mathbb{R}} \tilde{I}_{\infty}(u_n,h_n)   \|\rightarrow 0\ \text{as}\ n\rightarrow \infty$,\\
(3) $dist((u_n,h_n),D_n\times\{1\})\rightarrow 0\ \text{as}\ n \rightarrow \infty.\\                $
Then $v_n :=h_n \star u_n \in S_c$ is a Palais-Smale sequence for $I_\infty$ on $S_c$ at $m_c$, and there exist Lagrange multipliers $\lambda_n \in \mathbb{R}$ such that
$$ I_{\infty}(v_n)\rightarrow m_c,\ \ \ P_{\infty}(v_n)\rightarrow 0, $$
$$ \| I_{\infty}'(v_n)+\lambda_n G'(v_n) \|_{(H^1(\mathbb{R}^N))^*} \rightarrow 0,\ \ \ \text{where}\ G(u)=\frac{1}{2} \int_{\mathbb{R}^N}u^2 dx,                     $$
as $n \rightarrow \infty$. So, combining those properties we can infer that
$$\frac{N(p-2)-4}{2N(p-2)} a \| \nabla v_n \|_2^2 + \frac{N(p-2)-8}{4N(p-2)}b \| \nabla v_n \|_2^4  \rightarrow m_c >0,\ \text{as}\ n \rightarrow \infty,  $$
and
$$ \aligned -\lambda_n c^2 &= a\| \nabla v_n \|_2^2 +b\| \nabla v_n \|_2^4 -\| v_n \|_p^p \\
                         &=\frac{N(p-2)-2p}{2p}\| v_n \|_2^2 = \frac{N(p-2)-2p}{N(p-2)}(a \|\nabla v_n \|_2^2 + b\| \nabla v_n \|_2^4).  \endaligned $$
Therefore, $v_n$ is bounded in $H^1(\mathbb{R}^N)$ and $\lambda_n$ is bounded in $\mathbb{R}$. We may assume that $ v_n \rightharpoonup v$ in $H^1(\mathbb{R}^N)$, $\|\nabla v_n\|_2^2 \rightarrow A^2$, and $\lambda_n \rightarrow \lambda > 0$. In fact, $\{ v_n \}$ is a (PS) sequence for $I_{\infty,\lambda}$ at level $m_c+\frac{1}{2}c^2$.

As a consequence of Proposition 3.1, $v_n$ can be rewritten as
$$  v_n = v + \sum_{j=1}^{k} w^j(\cdot-y_n^j)+o(1)    $$
in $H^1(\mathbb{R}^N)$, where $k \geq 0$ and $w^j \neq 0$, $v$ are solutions to
\begin{equation}\label{equ limit}
  -(a+bA^2)\triangle w + \lambda w = |w|^{p-2}w
\end{equation}
and $|y_n^j|\rightarrow \infty$. Moreover, from (\ref{dec equ}) and (\ref{mp level des}), we get
\begin{equation}\label{A des identy}
 \aligned
& c^2 = \| v \|_2^2 + \sum_{j=1}^{k}\| w^j \|_2^2,\\
& A^2=\|\nabla v \|_2^2 + \sum_{j=1}^{k} \| \nabla w^j \|_2^2,\\
&I_{\infty,\lambda}(v_n) \rightarrow J_{\infty}(v) + \sum_{j=1}^{k}J_{\infty}(w^j),
\endaligned
\end{equation}
and hence,
$$ I_{\infty}(v_n)\rightarrow m_c =  J_{\infty}(v)-\frac{\lambda}{2}\|v \|_2^2 + \sum_{j=1}^{k}J_{\infty}(w^j)-\frac{\lambda}{2}\sum_{j=1}^{k}\| w^j \|_2^2.                $$

If $v \neq 0$ and $k\geq1$, combining $v$ and $w^j$ solves (\ref{equ limit}), we have the Pohozaev type identities
$$ \tilde{P}(v):= (a+bA^2)\| \nabla v \|_{2}^2 - \frac{N(p-2)}{2p}\| v \|_p^p=0,   $$
$$ \tilde{P}(w^j):= (a+bA^2)\| \nabla w^j \|_{2}^2 - \frac{N(p-2)}{2p}\| w^j \|_p^p,\ \ \text{for}\ j=1,2,\cdots,k.    $$
Since $v \neq 0$, $k\geq1$ and (\ref{A des identy}), we have $\| \nabla v \|_2^2 < A^2$, so that $P_{\infty}(v)<0$, that is,
$$ a\| \nabla v \|_{2}^2 + b\| \nabla v  \|_2^4  < \frac{N(p-2)}{2p}\| v \|_p^p.  $$
Then, there exist a constant $s_1 \geq 1$ such that $I_{\infty}(s_1 \star v)<0$ since $\frac{N(p-2)}{2p} \leq \frac{4}{p}$ and
$$  I_{\infty}(s\star v)=\frac{a s^2 }{2}\| \nabla v \|_2^2 +\frac{s^4}{4}\left[ b\| \nabla v \|_2^4 -\frac{4}{p}s^{\frac{N(p-2)}{2}-4}\|v\|_p^p  \right].     $$

We can define a path $\sigma^*(t):= ((1-t)s_2+ts_1)\star v \in \Sigma_{\| v \|_2}$ considering $\| \nabla (s_2 \star v) \|_{2}^2 < r(\| v\|_2)$ for some $0 < s_2 < 1$, and there exists a constant $t^* \in [0,1]$ such that $I_{\infty}(\sigma(t^*))=\max\limits_{t \in [0,1]}I_{\infty}(\sigma(t))$.
In particular, $P_{\infty}(\sigma(t^*))=0$ and $s_*:=((1-t^*)s_2+t^*s_1) < 1$ since $P_{\infty}(v)<0$.

\noindent\textbf{Claim:}\ If $p=2+\frac{8}{N}$, then $\| v \|_2:=\alpha > c_0$.\\
\noindent Define a auxiliary function $f(t):= a t + bt^2 - 2\alpha^{\frac{8}{N}-2}   \| Q_{2+\frac{8}{N}} \|_2^{-\frac{8}{N}} t^2 $, and
$$ \aligned P_{\infty}(v) &= a\| \nabla v \|_2^2 + b\| \nabla v \|_2^4- \frac{2N}{N+4} \|v \|_p^p  \\
                          &\geq a\| \nabla v \|_2^2 + b\| \nabla v \|_2^4- 2\alpha^{\frac{8}{N}-2}   \| Q_{2+\frac{8}{N}} \|_2^{-\frac{8}{N}} \|\nabla v \|_2^4  \\
                          &=f(\| \nabla v \|_2^2).
\endaligned $$
Using $P_{\infty}(v)<0$, we have $\min_{t \geq 0} f(t)<0$, so that $\alpha > \left(\frac{b}{2}\right)^{\frac{8}{8-2N}}\| Q_{2+\frac{8}{N}} \|_2^{\frac{8}{8-2N}}=c_0$.\\
Thus, $m_{\|v\|_2}$ is well defined. However, we have
\begin{equation}\label{the splitting value}
\aligned J_{\infty}(v)-\frac{\lambda}{2}\|v \|_2^2 &=J_{\infty}(v)-\frac{\lambda}{2}\|v \|_2^2 -\frac{2}{N(p-2)}\tilde{P}(v)\\
= &\frac{N(p-2)-4}{2N(p-2)}a \|\nabla v \|_2^2 + \frac{N(p-2)-8}{4N(p-2)}bA^2\|\nabla v\|_2^2   \\
> &\frac{N(p-2)-4}{2N(p-2)}a \|\nabla v \|_2^2 + \frac{N(p-2)-8}{4N(p-2)}b\|\nabla v\|_2^4      \\
> &s_*^2\frac{N(p-2)-4}{2N(p-2)}a \|\nabla v \|_2^2 + s_*^4 \frac{N(p-2)-8}{4N(p-2)}b\|\nabla v\|_2^4      \\
= & I_{\infty}(\sigma^*(t^*))
\geq m_{\| v \|_2}.
\endaligned \end{equation}

\noindent Similarly, we have $\| w^j \|_2 > c_0$ and $ J_{\infty}(w^j)-\frac{\lambda}{2}\|w^j \|_2^2 > m_{\|w^j \|}$ for all $1\leq j\leq k$.
Now, from \emph{Lemma 2.2(2)}, we get a contradiction
$$ m_c + o(1)\geq m_{\|v \|_2} +\sum_{j=1}^{k} m_{\| w^j \|} + o(1)\geq (k+1)m_c+o(1).$$

\noindent Therefore, either $k=1$ and $v=0$, or $k=0$ and $v \neq 0$.

If $k=1$ and $v=0$, then $v_n(\cdot+y_n^1)+o(1) = w^1$. On the other hand, due to point (3) that $dist((u_n,h_n),D_n\times\{0\}) \rightarrow 0$, we obtain
$$ \beta(w^1)=\beta (v_n(\cdot+y_n^1))+o(1)=y_n^1 +o(1)   $$
which contradicts the fact $\beta$ is continuous and $| y_n^1|\rightarrow \infty$.

If $k=0$ and $v \neq 0$, then $v_n \rightarrow v$ in $H^{1}(\mathbb{R}^N)$. Using again point (3), we also have $\beta(v)=0$. Hence, by the uniqueness, $v_n \rightarrow \pm u_c$ in $H^1(\mathbb{R}^N)$. This implies
$$ I(v_n)=I_{\infty}(v_n)+\frac{1}{2}\int_{\mathbb{R}^N}V(x)v_n^2 dx \rightarrow m_c + \frac{1}{2}\int_{\mathbb{R}^N}V(x)v_n^2 dx > m_c          $$
which is a contradiction. \qed

\begin{Lem}
For any $c > c^*$, then $m_{V,c} \geq L_c$ holds.
\end{Lem}
\noindent\textbf{Proof:} Similar to\cite[Proposition 3.5]{Bartsch2021} we omit it. \qed

\begin{Lem}
For any $c > c^*$ and for any $\varepsilon > 0$, there exist $\bar{R}$ and $\bar{h}>1$ such that for $Q=B_R \times [h_1,h_2]$ with $R\geq\bar{R}$ and $0<h_1 \leq \bar{h}^{-1} < \bar{h} \leq h_2$, the following holds:
$$    \max_{(y,h)\in \partial Q} I(h\star u_c(\cdot - y)) < m_c + \varepsilon.        $$
\end{Lem}
\noindent\textbf{Proof:} Similar to \cite[Proposition 3.6]{Bartsch2021} and \cite[Section 4]{Bartsch2021}, we omit it.\qed

By $Lemma\ 3.5$ and $Lemma\ 3.6$, we may choose $R>0$ and $0<h_1<h_2$ such that
$$  \max_{(y,h)\in\partial Q}I(h\star u_c(\cdot -y))< m_{V,c}. $$
This implies that $I$ has a linking geometry and there exists a Palais-Smale sequence at level $m_{V,c}$. In order to estimate $m_{V,c}$, we have the following Lemma.
\begin{Lem}
Assume $(V1)$ holds, for $c>c^*$, and $h_1>0$ small enough, $h_2>0$ large enough, then
$$   m_{V,c}\leq m_c + \frac{1}{2}\| V \|_{\infty}c^2. $$
\end{Lem}
\noindent\textbf{Proof:}
$$\aligned  m_{V,c} &\leq \max_{(y,h)\in Q} \left\{I_{\infty}(h\star u_c(\cdot-y))+\int_{\mathbb{R}^N}V(x)(t\star u_c)^2(x-y)dx\right\} \\
                    &\leq \max_{h\in [h_1,h_2]} I_{\infty}(h\star u_c) + \frac{1}{2}\|V \|_{\infty}c^2 \\
                    &=m_c+\frac{1}{2}\|V \|_{\infty}c^2\endaligned $$
provided $h_1>0$ small enough and $h_2>0$ large enough. \qed

\begin{Lem}
Assume $(V2)$ holds, for $c>c^*$, set $\bar{\nu} = a^{-1}S^{-2} \| V\|_{\frac{3}{2}}$ and $\nu= \frac{3(p-2)\bar{\nu}}{3p-10-4\bar{\nu}}$. Then $0<\nu<\frac{3}{2}$ and
$$m_{V,c} \leq (1+\nu)m_c. $$
\end{Lem}
\noindent\textbf{Proof:} We have

\begin{equation}\label{poles m}
\aligned m_{V,c} \leq  &\max_{h>0,y\in \mathbb{R}^3} I(h\star u_c(\cdot-y)) \\
= &\max_{h>0,y\in \mathbb{R}^3} \left[ I_{\infty}(h\star u_c) +\frac{1}{2}\int_{\mathbb{R}^3}V(x+y)(h\star u_c)^2 dx \right]  \\
\leq & \max_{h>0}\left[ I_{\infty}(h\star u_c)+\frac{1}{2}\|V \|_{\frac{3}{2}}\|h\star u_c \|_{6}^2  \right]  \\
\leq &\max_{h>0}\left[ I_{\infty}(h\star u_c)+\frac{1}{2}a \bar{\nu} \|\nabla u_c \|_{2}^2 h^2 \right] \\
= & \max_{h>0}\left[ \frac{1}{2}a(1+\bar{\nu}) \|\nabla u_c \|_{2}^2 h^2 + \frac{1}{4} b h^4 \| \nabla u_c \|_2^4 -\frac{1}{p}\| u_c \|_p^p h^{\frac{3(p-2)}{2}}   \right]. \endaligned \end{equation}
It is easy to see $ \bar{\nu}<1 $ so that $\nu>0$ and $\nu<\frac{2}{3}$ provided (\ref{V21}). Define two auxiliary functions
$$ f(h) = \frac{a}{2}\|\nabla u_c \|_2^2t^2 + \frac{b}{4}\| \nabla u_c \|_2^4 t^4 -\frac{1}{p}\|u_c\|_p^p t^{\frac{3(p-2)}{2}} $$
and
$$ g(t) = (1+\bar{\nu})\frac{a}{2} \|\nabla u_c \|_2^2t^2+ \frac{b}{4}\| \nabla u_c \|_2^4 t^4 -\frac{1}{p}\|u_c\|_p^p t^{\frac{3(p-2)}{2}}$$
for all $t>0$. By the definition of $m_c$ and $u_c$, we deduce $\max_{h>0}f(t) = f(1) = m_c$ and $f'(t) \leq 0$ if $t \geq 1$. Since $\bar{\nu}>0$, then we must have $t>1$ if $h'(t)\leq 0$.

Consider $$(1+\nu)f(t)-g(t)= \nu\left[(1-\frac{\bar{\nu}}{\nu})\frac{a}{2}\| \nabla u_c \|_2^2 +\frac{b}{4}\| \nabla u_c \|_2^4 -\frac{1}{p} \| u_c \|_p^p t^{\frac{3(p-2)}{2}}\right],$$
then there exist a unique $t_1 >0$ such that $(1+\nu)f(t_1)-g(t_1)$ and $(1+\nu)f(t)-g(t)>0$ if $t>0$ small enough since $\bar{\nu} < \nu$. Moreover, we have
$$\aligned  g'(t_1)t_1 & =(1+\bar{\nu})a \| \nabla u_c \|_2^2 t_1^2  + b \| \nabla u_c \|_2^4 t_1^4 -\frac{3(p-2)}{2p}\| u_c \|_p^p t_1^\frac{3(p-2)}{2} \\
                       & = a\| \nabla u_c \|_2^2 t_1^2\left( \frac{4-3(p-2)}{4}+ \frac{\bar{\nu}}{\nu}\frac{3(p-2)}{4}+\bar{\nu} \right) +b\|\nabla u_2 \|_2^4t_1^4\left(\frac{8-3(p-2)}{8}\right) \\
                       & \leq 0,\endaligned $$
and so that $t_1 \geq 1 $ and $f'(t_1) \leq 0 $. Thus, we deduce $\max_{t>0}g(t) \leq \max_{t>0}(1+\nu)f(t)=(1+\nu)m_c$ which
completes the proof. \qed

\subsection{Under radial symmetry assumption}
Furthermore, we consider

\noindent\textbf{Case 1:}  $p \in [2+\frac{8}{N},2^*)$, $c>c^*$, and $V(x)$ satisfies $(V4)$, $ 2 \leq N \leq 3$ or $(V4')$, $N=3$.

\noindent\textbf{Case 2:} $V(x)$ satisfies $(V5)$, $N\geq 2$, $2+\frac{4}{N} < p <\min\{2+\frac{8}{N},2^* \}$, and $c>c_1$.

\noindent Namely, let $ U_c$ be $u_c$ of \emph{Lemma\ 2.2} in \textbf{Case 1} and be $u_{c,1}$ of \emph{Lemma\ 2.3} in \textbf{Case 2}.

\noindent We define
$$ S_{c}^r =\{ u \in H_{rad}^1(\mathbb{R}^N):\| u \|_2^2=c^2 \}, $$
$$ \Psi_c =\{  \psi \in \mathcal{C}([0,1],S_c^r):\phi(0)=h_1 \star U_c, \psi(1) = h_2 \star U_c   \}          $$
where $0<h_1<h_2$ such that
$$  \max\{  I(h_1 \star U_c), I(h_2 \star U_c)   \}< M_{c}:=\left\{ \aligned  &m_c,\ \ \ \ \ \text{in}\ \textbf{Case 1},   \\
                                                                          &m_{c,2},\ \ \ \text{in}\ \textbf{Case 2}.              \endaligned  \right. $$
The choice of $h_1$ and $h_2$ is easy in \textbf{Case 1}, because $I$ has mountain pass structure and
$$ I(h\star U_c) \rightarrow \left\{ \aligned    &0,\ \ \ \ \ \ \ \text{as}\ h \rightarrow 0    \\
                                                 &-\infty,\ \ \text{as}\ h \rightarrow \infty.   \endaligned     \right.         $$
In \textbf{Case 2}, we can choose $h_1$ small enough and $h_2 = 1$, since we have
$$ I(U_c)=m_{c,1}+\frac{1}{2}\int_{\mathbb{R}_N}V(x)U_c^2 dx < m_{c,1} + \frac{1}{2}\| V \|_{\infty}c^2 < M_c $$
from $ \| V\|_{\infty} <2c^{-2}(m_{c,2}-m_{c,1})$ (by $(V5)$).

We look for a critical point of $I$ at the level
$$ M_{V,c}: = \inf_{\psi \in \Psi}\max_{t\in[0,1]}I(\psi(t)). $$

\begin{Lem}
In both \textbf{\emph{Case 1}} and \textbf{\emph{Case 2}}, we have $M_{V,c}>m_c$.
\end{Lem}
\noindent\textbf{Proof:} Using the fact $I_{\infty}(u^*) \leq I_{\infty}(u)$ while $u^*$ is the radially decreasing rearrangement of $u \in H^{1}(\mathbb{R}^N)$, we have
$$ M_c=\bar{M}_c:= \inf_{\phi \in \Psi}\max_{t \in [0,1]} I_{\infty}(\phi(t)).    $$
Clearly, $M_{V,c} \geq M_c$ since $V \geq 0$. If the equality holds, then there would exist a sequence $\psi_n \in \Psi$ such that
$$    0 \leq \max\limits_{t\in[0,1]}I(\psi(t))-M_c \leq \frac{1}{n}.         $$
Arguing as in the proof of Lemma 3.4, it would be possible to construct a Palais-Smale sequence $v_n \in H_{rad}^1(\mathbb{R}^N)$ such that $v_n \rightarrow \pm U_c$ strongly in $H^1(\mathbb{R}^N)$ and
$$  dist(v_n,\psi_n([0,1]))=\| v_n - \psi_n (\bar{t}_n)    \|_{H^1(\mathbb{R}^N)} \rightarrow 0,    $$
for some $\bar{t}_n \in [0,1]$. Therefore,
$$ \aligned \max\limits_{t\in[0,1]}I(\psi_n(t)) \geq I(\psi_n(\bar{t}_n)) &= F(v_n) + o(1) \\
                                                                          &=I_{\infty}(v_n)+\frac{1}{2}\int_{\mathbb{R}^N} V(x)v_n^2dx +o(1)\\
                                                                          &\rightarrow m_c +\frac{1}{2}\int_{\mathbb{R}^N} V(x)v_n^2dx\\
                                                                          &>m_c \endaligned $$
a contradiction. \qed

\section{Proof of Theorem 1.2, 1.3 and 1.4}
In this section, we always assume that $1 \leq N \leq 3$, $c>c_*$ and $p\in[2+\frac{8}{N},2^*)$ hold. The proof could be divided into two steps. Firstly, we construct a bounded Palais-Smale sequence of $I$ at $m_{V,c}$ by adopting the approach from \cite{Jeanjean1997} and \emph{Lemma 3.3}. Secondly, we show that the sequence converge strongly in $H^1(\mathbb{R}^N)$ drawing support from \emph{Proposition 3.1}.

We define a auxiliary $\mathcal{C}^1$ functional
$$\tilde{I}(u,h):=I(h\star u) \ \ \text{for all}\ (u,h)\in H^1(\mathbb{R}^N)\times\mathbb{R},$$
$$ \tilde{\Phi}_c:=\left\{ \tilde{\phi}\in\mathcal{C}(Q,S_c\times \mathbb{R}): \tilde{\phi}(y,h):=(h\star u_c(\cdot-y),1)\ \text{for all}\ (y,h)\in\partial Q  \right\},  $$
and
$$ \tilde{m}_{V,c}:= \inf_{\tilde{\phi}\in\tilde{\Phi}_c}\max_{(y,h)\in Q}\tilde{I}(\tilde{\phi}(y,h)). $$

\begin{Lem}
\noindent(1)\ $\tilde{m}_{V,c}=m_{V,c}$.\\
\noindent(2)\ If $(u_n, h_n)$ is a $(PS)$ sequence for $\tilde{I}$ at level $m \in\mathbb{R}$ and $h_n \rightarrow 1$, then $(h_n \star u_n)_n$ is a $(PS)$ sequence for $I$ at level $m$.
\end{Lem}
\noindent\textbf{Proof:}\ (1)\ Since $\Phi_c \times \{ 1 \} \subset \tilde{\Phi}_c$, then $\tilde{m}_{V,c} \geq m_{V,c}$. On the other hand, for any $\tilde{\phi}=(u,h) \in  \tilde{{\Phi}}_c$, the function $\phi:=h \star u \in \Phi_c$ satisfies
$$   \max_{(y,h)\in Q}\tilde{I}(\tilde{\phi}(y,h))=\max_{(y,h)\in Q}I(\phi(y,h))              $$
so that $m_{V,c} \leq \tilde{m}_{V,c}$.

(2) The proof is similar to that of \cite{Jeanjean1997} and is omitted. \qed

\begin{Prop}
Let $\tilde{\phi}_n \in \tilde{\Phi}_c$ be a sequence such that
$$  \max_{(y,h)\in Q}\tilde{I}(\tilde{\phi}_n(y,h)) \leq m_{V,c}+\frac{1}{n}.    $$
Then, there exist a sequence $(u_n,h_n)\in S_c \times \mathbb{R}$ and $d>0$ such that
$$  m_{V,c}- \frac{1}{n} \leq \tilde{I}(u_n,h_n)   \leq m_{V,c}+\frac{1}{n}      $$
$$  \min_{(y,h)\in Q}\| (u_n,h_n)-\tilde{\phi}_n(y,h)  \|_{H^1(\mathbb{R}^N) \times \mathbb{R}}  \leq \frac{d}{\sqrt{n}}  $$
$$  \| \nabla_{S_c \times \mathbb{R}} \tilde{I}(u_n,h_n)     \| \leq \frac{d}{\sqrt{n}}.      $$
The last inequality means:
$$ |D\tilde{I}(u_n,h_n)[(z,s)]|\leq \frac{d}{\sqrt{n}} (\|z\|_{H^1(\mathbb{R}^N)} + |s|)    $$
for all
$$ (z,s) \in \left\{ (z,s) \in H^1(\mathbb{R}^N) \times \mathbb{R} : \int_{\mathbb{R}^N}z u_n dx = 0 \right\}.   $$
\end{Prop}

\noindent\textbf{Proof:}\ Apply \emph{Lemma 3.3} to $\tilde{I}$ with
$$ E:=S_c \times \mathbb{R},\ \ \ T:=\{ (h \star u_c (\cdot - y),0):(y,h)\in \partial Q  \},  $$
$\ \ \ \ \ \ \ \ \ \ \ \ \ \ \ \ \ \ \ \ \ \ \ \ \ \ \ \mathcal{X}:=\tilde{\Phi}_c,\ \ X_n := \{\tilde{\phi}_n(y,h):(y,h)\in Q \}. $\qed

As a consequence we obtain a bounded Palais-Smale sequence for $I$ at the level $m_{V,c}$.
\begin{Prop}
Assume that $(V1)$ hold, then there exists a bounded sequence $(v_n)_n$ in $S_c$ such that
\begin{equation}\label{equ:PS}
  I(v_n)\rightarrow m_{V,c},\ \ \ \nabla_{S_c}I(v_n)\rightarrow 0
  \end{equation}
and
\begin{equation}\label{equ:P condition}
  a\| \nabla v_n \|_2^2 + b \| \nabla v_n \|_2^4 - \frac{N(p-2)}{2p}\| v_n \|_p^p - \frac{1}{2}\int_{\mathbb{R}^N} V(x)(N v_n^2 + 2v_n \nabla v_n \cdot x)dx  \rightarrow 0
\end{equation}
as $n \rightarrow \infty$. Moreover, the sequence of Lagrange multipliers satisfies, up to subsequence,
\begin{equation}\label{equ:lambda V}
  \lambda_n:= -\frac{DI(v_n)[v_n]}{c^2} \rightarrow \lambda>0.
\end{equation}
\end{Prop}
\noindent\textbf{Proof:}\ First we choose a sequence $\phi_n \in \Phi_c$ such that
$$ \max_{(y,h)\in Q}I(\phi_n(y,h)) \leq m_{V,c}+\frac{1}{n}.  $$
Since $I$ is even, we can assume that $\phi_n(y,h)\geq 0$ almost everywhere in $\mathbb{R}^N$. Applying \emph{Proposition 4.2} to $\tilde{\phi}_n(y,h):=(\phi_n(y,h),1)\in \tilde{\Phi}_c$, we can prove the existence of a sequence $(u_n,h_n) \in H^1(\mathbb{R}^N)\times \mathbb{R}$ such that $I(h_n \star u_n) \rightarrow m_{V,c}$. Note that, by \emph{Proposition 4.2}, we have
$$  \min_{(y,h)\in Q}\|  (u_n,h_n)-\tilde{\phi}_n(y,h)  \|_{H^1(\mathbb{R}^N)\times \mathbb{R}}    \leq \frac{d}{\sqrt{n}}   $$
so that $h_n \rightarrow 1$ as $n \rightarrow \infty$ and there exists $(y_n, h_n^*) \in B_R \times [h_1,h_2]$ such that $u_n - \phi_n(y_n,h_n^*)=o(1)$ in $H^1(\mathbb{R}^N)$, as $n \rightarrow \infty$.

We define $$v_n:=h_n \star u_n.$$
Observe that, since $\phi_n(y_n,h_n^*)\geq 0$ a.e. in $\mathbb{R}^N$, then $\| u_n^- \|_2 \leq \| u_n -\phi_n(y_n,h_n^*) \|_2 =o(1)   $
and we can deduce $u_n^- \rightarrow 0$ a.e., up to a subsequence. So
\begin{equation}\label{positive ps}
  \|v_n^{-1} \|_2 \rightarrow 0,\ \ \ \text{ as } n \rightarrow \infty.   \end{equation}
Moreover, by \emph{Lemma 4.1} $(v_n)_n$ is a Palais-Smale sequence for $I$. Similarly we have for $w \in H^1(\mathbb{R}^N)$, setting $\tilde{w}_n:=(h_n^{-1}\star w)$, such that
$$  \nabla (I-I_{\infty})(v_n)[w]=\int_{\mathbb{R}^N} V(h_n^{-1}x)u_n \tilde{w}_n dx, $$
which implies
$$    DI(v_n)[w]=D\tilde{I}(u_n,h_n)[(\tilde{w}_n,0)]+o(1)\| \tilde{w}_n \|.               $$
Moreover, it is easy to see that $\int_{\mathbb{R}^N}v_nw =0$ is equivalent to $\int_{\mathbb{R}^N}u_n\tilde{w}=0$. Since $\| \tilde{w}_n\|_{H^1(\mathbb{R}^N)}^2 \leq 2\| w \|_{H^1(\mathbb{R}^N)}^2$ for $n$ large, we have (\ref{equ:PS}).

Using the fact that $D\tilde{I}(u_n,h_n)[(0,1)]\rightarrow 0$ and $ h_n \rightarrow 1 $ as $n\rightarrow \infty$, we deduce that
$$ \partial_{h}\left( \int_{\mathbb{R}^N}V(x)h^{N}u^2(hx) dx \right) =\int_{\mathbb{R}^N}V(x)(Nh^{N-1}u^2(hx)+2h^{N}u^2(hx)\nabla u(hx)\cdot hx)dx  $$
and
$$  \partial_h I_{\infty}(h \star u) = a \| \nabla v \|_2^2 + b \| \nabla v \|_2^4 - \frac{N(p-2)}{2p}\| v \|_p^p. $$
So that, (\ref{equ:P condition}) holds.

Now we prove $(v_n)_n$ is bounded in $H^1(\mathbb{R}^N)$. Recall $\theta = N(p-2)-4$, $\eta=N(p-2)-8$ and set
$$   A_n=\| \nabla v_n \|_2^2,\ \ B_n=\| v_n \|_p^p,\ \ C_n=\int_{\mathbb{R}^N}V(x)v_n^2 dx,\ \ D_n=\int_{\mathbb{R}^N}V(x)v_n(\nabla v_n \cdot x)dx.     $$
Then (\ref{equ:PS}) (\ref{equ:P condition}) and (\ref{equ:lambda V}) can be expressed in the form
$$  a A_n +\frac{b}{2}A_n^2+ C_n -\frac{2}{p}B_n = 2m_{V,c} + o(1) \ \ \ \ \text{as }n\rightarrow \infty            $$
$$  aA_n + bA_n^2 -\frac{N(p-2)}{2p}B_n +\frac{N}{2}C_n + D_n = o(1)  \ \ \ \ \text{as }n\rightarrow \infty           $$
$$  aA_n + bA_n^2 + \lambda_n c^2 + C_n = B_n +o(1)(A_n^{1/2}+1)  \ \ \ \ \text{as }n\rightarrow \infty.            $$
We obtain
$$ \frac{N(p-2)-4}{2p}B_n = 2m_{V,c}+D_n+\frac{N-2}{2}C_n+\frac{b}{2}A_n^2 + o(1).   $$
As a consequence
$$   aA_n + \frac{b}{2}\frac{\eta}{\theta} A_n^2 = 2\frac{N(p-2)}{\theta}m_{V,c}+\frac{4}{\theta}D_n +\frac{N(4-p)}{\theta}C_n + o(1).    $$
This implies
$$  \aligned  a\theta A_n + \frac{b}{2}\eta A_n^2 & = 2N(p-2)m_{V,c}+4D_n+N(4-p)C_n+o(1)\\
                                                 & \leq 4N(p-2)m_c + 4\| W \|_{\infty}cA_n^{1/2}  \endaligned$$
using the fact that $p>4$, $m_{V,c} < 2m_c$, $C_n \geq 0$, $W(x)=V(x)|x|$ and H\"{o}lder inequality, so that $A_n$ is bounded.
More precisely, we get \\
\begin{equation}\label{A bound} A_n^{1/2} \leq \frac{2\|W\|_{\infty}c + 2\sqrt{\| W \|_{\infty}^2 c^2 + a\theta N(p-2)m_c       }         }{a\theta}. \end{equation}

To conclude, we prove $\lambda_n$ admits a subsequence that converges to a positive bounded limit.
Since $A_n$ and $C_n$, then $D_n$ and $B_n$ is bounded due to H\"{o}lder inequality, therefore also $\lambda_n$ is bounded. Up to a subsequence, we may assume
$$ A_n \rightarrow A \geq 0,\ \ B_n \rightarrow B \geq 0,\ \ C_n \rightarrow C \geq 0,\ \ D_n \rightarrow D \in \mathbb{R},\ \ \lambda_n \rightarrow \lambda \in \mathbb{R}.         $$
Passing to the limit, correspondingly, we have
$$  a A +\frac{b}{2}A^2+ C -\frac{2}{p}B = 2m_{V,c}   $$
$$  aA + bA^2 -\frac{N(p-2)}{2p}B +\frac{N}{2}C + D = 0           $$
$$  aA + bA^2 + \lambda c^2 +C = B         $$
which imply
$$ \aligned \lambda c^2 & =B-aA-C-bA^2 \\
                        & =\frac{p-2}{p}B-2m_{V,c}-\frac{b}{2}A^2 \\
                        & =\frac{p(2-N)+2N}{N(p-2)-4}2m_{V,c} + \frac{2(p-2)}{N(p-2)-4}D + \frac{(p-2)(N-2)}{N(p-2)-4}C + \frac{p(2-N)+2N}{N(p-2)-4}\frac{b}{2}A^2.
\endaligned $$
Notice that
$$ \frac{(p-2)(N-2)}{N(p-2)-4}\geq 0, \ \ \text{when }N \geq 2, \ \   \frac{(p-2)(N-2)}{N(p-2)-4}=\frac{2-p}{p-6}<0,\ \ \text{when }N=1.    $$
so that $\lambda > 0$ provided
$$ (p-2)|D| \leq \left\{ \aligned &  4 m_c,\ N=1,   \\
                                   &\bigl[ p(2-N)+2N \bigr]  m_c,\ N=2,3.   \endaligned  \right.$$
Using H\"{o}lder inequality it is possible to see that
$$  (p-2)|D| \leq  (p-2) \|W\|_{\infty}c A^{\frac{1}{2}} \leq \left\{  \aligned &  4 m_c,\ N=1,   \\
                                                                                  &\bigl[ p(2-N)+2N \bigr]  m_c,\ N=2,3.   \endaligned      \right.   $$
provided
\begin{equation}\label{W bound}
 \aligned & \| W \|_{\infty} \leq \frac{m_c^{1/2}}{c}\left(\frac{4a(p-6)}{(p-2)^3+4(p-2)}\right)^{1/2}, \ \ N=1,\\
          & \| W \|_{\infty} \leq \frac{m_c^{1/2}}{c}\left(  \frac{a(N(p-2)-4)[p(2-N)+2N]^2}{4(p-2)(N(p-2)+(p(2-N)+2N))}   \right)^{1/2}, \ \ N=2,3. \endaligned
\end{equation}
%
From the fact that $4<p(2-N)+2N$, $C \leq \| V \|_{\infty}c^2$ and $m_{V,c} \leq m_c +\frac{1}{2}\|V \|_{\infty}c^2$, we deduce
\begin{equation}\label{lambda bound}
\aligned  \lambda c^2\leq &\frac{p(2-N)+2N}{N(p-2)-4}2m_{V,c} + \frac{2(p-2)}{N(p-2)-4}D + \frac{(p-2)(N-2)}{N(p-2)-4}C + \frac{p(2-N)+2N}{N(p-2)-4}\frac{b}{2}A^2\\
                     \leq & 2q+\frac{8\|V \|_{\infty}}{N(p-2)-4}m_c + \frac{p(2-N)+2N}{N(p-2)-4}\frac{b}{2}A^2.
                             \endaligned
\end{equation}

\begin{Lem}
Assume that (V3) hold, for any $v \in S_{c}$ if $p>2+\frac{8}{N}$, and $v \in S_{c}$ satisfying $P_{\infty}(v) \leq 0$ if $p=2+\frac{8}{N}$, we have
$$I(v)-\frac{2}{N(p-2)}P(v) \geq m_{c}-\frac{2}{N(p-2)}t^{\frac{N(p-2)}{2}}P(v), \ \forall t > 0,$$
where
$$\aligned P(v)&=a\| \nabla v \|_2^2 + b \| \nabla v \|_2^4 - \frac{N(p-2)}{2p}\| v \|_p^p - \frac{1}{2}\int_{\mathbb{R}^N} V(x)(N v^2 + 2v \nabla v \cdot x)dx \\
               &=a\| \nabla v \|_2^2 + b \| \nabla v \|_2^4 - \frac{N(p-2)}{2p}\| v \|_p^p - \frac{1}{2}\int_{\mathbb{R}^N} \langle\nabla V(x)\cdot x\rangle v^2dx. \endaligned  $$
\end{Lem}
\noindent\textbf{Proof:} A direct computation shows that
$$\aligned I(v)-I(t\star v)
=&\frac{1}{2} \int_{\mathbb{R}^N}\biggl[  V(x)-V\left(\frac{x}{t}\right)+\frac{2}{N(p-2)}\left(1-t^{\frac{N(p-2)}{2}}\right)\langle \nabla V(x)\cdot x \rangle  \biggr]v^2 dx\\ &+\left(\frac{1-t^2}{2}-\frac{2(1-t^{\frac{N(p-2)}{2}})}{N(p-2)}\right)a\|\nabla v\|_2^2+\left(\frac{1-t^4}{4}-\frac{2(1-t^{\frac{N(p-2)}{2}})}{N(p-2)} \right)b\|\nabla v\|_2^4 \\
& + \frac{2(1-t^{\frac{N(p-2)}{2}})}{N(p-2)}P(v).\endaligned$$
By standard argument and $(V3)$, we have, for all $t > 0$,
$$ V(x)-V\left(\frac{x}{t}\right)+\frac{2}{N(p-2)}\left(1-t^{\frac{N(p-2)}{2}}\right)\langle \nabla V(x)\cdot x \rangle \geq 0,$$
$$ \frac{1-t^2}{2}-\frac{2(1-t^{\frac{N(p-2)}{2}})}{N(p-2)} \geq 0,\ \  \frac{1-t^4}{4}-\frac{2(1-t^{\frac{N(p-2)}{2}})}{N(p-2)} \geq 0.$$
Then, $$ I(v)-\frac{2}{N(p-2)}P(v) \geq  I(t\star v)-\frac{2}{N(p-2)}t^{\frac{N(p-2)}{2}}P(v),\ \forall t>0. $$
If $P_{\infty}(v) \leq 0$, similar to the argument of \emph{Lemma 3.4}, we can construct a path $\sigma \in \Sigma_{c}$ in \emph{Lemma 2.2}, by the definition of $m_c$, so that
$$ I(v)-\frac{2}{N(p-2)}P(v) \geq  m_{c}-\frac{2}{N(p-2)}t^{\frac{N(p-2)}{2}}P(v),\ \  \forall t >0. $$\qed

\begin{Lem}
Assume $(V1)$ and $(V3)$ hold, then the sequence $(v_n)$ converges strongly to $v \in S_c$, after passing to a subsequence.
\end{Lem}
\noindent\textbf{Proof:}
Since $v_n$ is bounded, after passing to a subsequence it converges weakly in $H^1(\mathbb{R}^N)$ to $ v \in H^1(\mathbb{R}^N)$. By (\ref{positive ps}) and weak convergence, $v$ is a nonnegative weak solution of
\begin{equation}\label{equation v}
   -(a+ bA^2 )\triangle v + (\lambda+V(x))v=|v|^{p-2}v
\end{equation}
such that $\| v \|_2 \leq c$ where $A^2:=\lim\limits_{n \rightarrow \infty}\| \nabla v_n \|_2^2$. We note that $v_n$ is a bounded Palais-Smale sequence of $I_\lambda$ at level $m_{V,c} + \frac{\lambda}{2}c^2$, therefore, by \emph{Proposition 3.1}, there exists an integer $k \geq 0$, $k$ non-trivial solutions $ w^1, w^2, \dots, w^k $ to the equation
\begin{equation}\label{equ:A equation 2}
  -(a+bA^2)\triangle w + \lambda w = |w|^{p-2}w
\end{equation}
and $k$ sequences $\bigr\{y_{n}^j\bigr\}\in H^{1}(\mathbb{R}^N)$, $1 \leq j \leq k$, such that $|y_n^j|\rightarrow \infty$ as $n \rightarrow \infty$.

\noindent Moreover, we have
\begin{equation}\label{dec equ 2}
\begin{split}
    & v_n -\sum_{j=1}^{k}w^{j}(\cdot - y_{n}^{j} )\rightarrow v \ \ \text{in}\ \ H^{1}(\mathbb{R}^N), \\
    & \|  v_n \|_2^2 \rightarrow \| v \|_2^2 + \sum_{j=1}^{k}\| w^j \|_2^2, \ \ \  A^2=\|\nabla v \|_2^2 + \sum_{j=1}^{k} \| \nabla w^j \|_2^2,
\end{split}
\end{equation}
and
\begin{equation}\label{mp level des 2}
  I_{\lambda}(v_n) \rightarrow J_{\lambda}(v) + \sum_{j=1}^{k}J_{\infty}(w^j)
\end{equation}
as $ n \rightarrow \infty$. It remains to show $k=0$, so that $v_n \rightarrow v$ strongly in $H^1(\mathbb{R}^N)$ and we are done.
Thus, by contradiction, we can assume that $k\geq 1$, or equivalently $\gamma:= \| v \|_2 < c$.

First we exclude the case $v=0$. In fact, if $v=0$ and $k = 1$, we have $w^1>0$ and $\|w^1\|_2 = c$ and $\|\nabla w^1 \|_2^2 = A^2$ so that $(\ref{mp level des 2})$ would give $m_{V,c}=m_c$, which is not possible due to \emph{Lemma 3.5}.
On the other hand, if $k \geq 2$, similar to \emph{Lemma 3.4}, we get $I_{\infty}(w^j) \geq m_{\alpha_j}$$(\alpha_j := \|w^j \|_2)$ and $m_{\alpha_j} > m_c$. However, the condition $I(v_n)\rightarrow m_{V,c}$ would give $$ 2m_c \leq m_{V,c}\leq m_c +\frac{\| V \|_{\infty}c^2}{2}, $$
which contradicts assumption $(V1)$.

Therefore from now on we will assume $v \neq 0$. From (\ref{dec equ 2}) and $I(v_n) \rightarrow m_{V,c}$, we deduce
$$ m_{V,c}+\frac{\lambda}{2}c^2 = J_{\lambda}(v) + \sum_{j=1}^{k}J_{\infty}(w^j). $$
Using (\ref{dec equ 2}) and adopting the argument of (\ref{the splitting value}), we have
$$ m_{V,c} \geq J_{\lambda}(v)-\frac{\lambda}{2}\gamma^2 + \sum_{j=1}^{k}m_{\alpha_{j}}.$$

\noindent \textbf{Claim:} $J_{\lambda}(v)-\frac{\lambda}{2}\gamma^2 \geq m_{\gamma}$.
Since $v$ solves (\ref{equation v}), we have
$$\aligned \tilde{P}_{V}(v)&=a\| \nabla v \|_2^2 + bA^2 \| \nabla v \|_2^2 - \frac{N(p-2)}{2p}\| v \|_p^p + \frac{1}{2}\int_{\mathbb{R}^N} V(x)(N v^2 + 2v \nabla v \cdot x)dx \\
               &=a\| \nabla v \|_2^2 + bA^2 \| \nabla v \|_2^2 - \frac{N(p-2)}{2p}\| v \|_p^p - \frac{1}{2}\int_{\mathbb{R}^N} \langle\nabla V(x)\cdot x\rangle v^2dx \\
               &=0.\endaligned  $$
Then, from $A^2 > \| \nabla v \|_2^2$ and $(V2)$, we have $$P(v)\leq\tilde{P}_{V}(v)=0,$$
and $P_{\infty}(v)\leq 0$ if $p=2+\frac{8}{N}$,
so that, combining with \emph{Lemma 4.4}, we deduce
$$ \aligned &J_{\lambda}(v)-\frac{\lambda}{2}\gamma^2\\
= &J_{\lambda}(v)-\frac{\lambda}{2}\gamma^2 - \frac{2}{N(p-2)}\tilde{P}_{V}(v)\\
= &\left(\frac{N(p-2)-4}{2N(p-2)}a+ \frac{N(p-2)-8}{4N(p-2)}bA^2\right) \|\nabla v \|_2^2 +\frac{1}{2}\int_{\mathbb{R}^N}\left(V(x)+\frac{2}{N(p-2)}\langle\nabla V(x)\cdot x\rangle\right)v^2dx\\
\geq & I(v)-\frac{2}{N(p-2)}P(v)\\
\geq & m_{\gamma}-\frac{2}{N(p-2)}t^{\frac{N(p-2)}{2}}P(v) \geq m_{\gamma}.   \endaligned $$
Thus, we have
$$\aligned 2m_c > m_c +\frac{\| V \|_{\infty}c^2}{2} &> m_{V,c} \geq J_{\lambda}(v)-\frac{\lambda}{2}\gamma^2 + \sum_{j=1}^{k}m_{\alpha_{j}}\\
                      & \geq m_{\gamma} + \sum_{j=1}^{k}m_{\alpha_{j}} \geq (k+1)m_c,   \endaligned$$
which shows $k=0$ and $v_n$ converges strongly to $v$ in $H^1(\mathbb{R}^N)$, otherwise it contradicts \emph{Lemma 2.2 (2)}. \qed

\noindent\textbf{Proof of Theorem 1.4:}\ Combining (\ref{equ:lambda V}) and \emph{Lemma 4.5}, we obtain $(v,\lambda) \in H^1(\mathbb{R}^N)\times \mathbb{R}^+$ is a solution to $(K_{V,c})$. \qed

\begin{Lem}
Assume that $(V1)$ hold, then there exist a constant $L>0$ such that the sequence $(v_n)$ converges strongly to $v \in S_c$, after passing to a subsequence, if $bm_c< L$.

\end{Lem}
\noindent\textbf{Proof:}\ Similar to \emph{Lemma 4.5}, we can exclude $v=0$, and assume $k\geq 1$.
Using (\ref{dec equ 2}) and adopting the argument of (\ref{the splitting value}), we have
$$ m_{V,c} \geq J_{\lambda}(v)-\frac{\lambda}{2}\gamma^2 + \sum_{j=1}^{k}m_{\alpha_{j}} \geq J_{\lambda}(v)-\frac{\lambda}{2}\gamma^2 + m_{\alpha} $$
where $\alpha = \max_{j}\alpha_j$.
Let $\beta>0$ be such that $(u_{\beta}, \lambda_{\beta})$ is the unique solution, $\lambda = \lambda_{\beta}$ and $I_{\infty}(u_{\beta})=m_{\beta}$ according to \emph{Lemma 2.1} and \emph{Lemma 2.2 (3)}.
By the uniqueness of solution to ($\ref{dec equ 2})$ and a suitable scaling, we get
\begin{equation}\label{alpha and lambda}
 \aligned & \lambda=\frac{2N-p(N-2)}{N(p-2)\alpha^2}\| \nabla w\|_2^2(a+bA^2),\\
   &\alpha^{\frac{4N-2p(N-2)}{4-N(p-2)}}=\|\nabla w\|_2^2\left[ \frac{4(a+bA^2)}{N(p-2)} \right]^{\frac{4}{4-N(p-2)}}\|Q_p\|_2^{\frac{2(p-2)}{4-N(p-2)}}. \endaligned
\end{equation}
Using (\ref{equ:Root equation}), (\ref{equ: lambda c}) and (\ref{alpha and lambda}), we have
$$ D_\alpha^2(a+bD_{\alpha}^2)^{\frac{4}{4-N(p-2)}}= \| \nabla w\|_2^2(a+bA^2)^{\frac{4}{4-N(p-2)}}  $$
and
$$\lambda_\alpha = \frac{2N-p(N-2)}{N(p-2)\alpha^2}D_{\alpha}^2(a+bD_{\alpha}^2)  $$
where $D_{\alpha}=\| \nabla u_\alpha \|_2$.
Considering a function $f(t)=t^2(a+bt^2)^{\frac{4}{4-N(p-2)}}-\| \nabla w\|_2^2(a+bA^2)^{\frac{4}{4-N(p-2)}}$, it is easy to see that $f(t)$ is non-increasing, $f(\|\nabla w \|_2)>0$ and $f(A)>0$. So $D_\alpha^2 \leq \| \nabla w \|_2^2 \leq A^2$
and
$$\frac{\lambda_{\beta}}{\lambda_{\alpha}} = \frac{\| \nabla w \|_2^2}{D_{\alpha}^2}\frac{(a+bA^2)}{(a+bD_{\alpha}^2)} \geq 1.$$ Thus, we have $\lambda_{\alpha} \leq \lambda_{\beta} $ and $\beta \leq \alpha $ follows.

Moreover,
$$ \aligned I_{\infty,{\lambda}}(u_\beta)&= \inf_{u \in H^1(\mathbb{R}^N)\setminus \{0\}} \max_{t \geq 0} I_{\infty,{\lambda}}(t u) \\
                                        &=\inf\{ I_{\infty,{\lambda}}(u): N_{\infty,\lambda}(u)=0  \}.  \endaligned  $$
Since $v$ solves (\ref{equation v}), we deduce
$$ \tilde{N}(v)=(a+bA^2)\|\nabla v \|_2^2 + \int_{\mathbb{R}^N}(V(x)+\lambda)v^2dx -\|v \|_p^p=0, $$
Thus, we have $N_{\infty,\lambda}(v)\leq \tilde{N}(v) =0$ by $V(x)\geq 0$ and $A^2 \geq \|\nabla v\|_2^2$, and there exist $t_* \leq 1$ such that $N_{\infty,\lambda}(t_*v)=0$.
Moreover,
\begin{equation}
\aligned J_{\lambda}(v)&=J_{\lambda}(v)-\frac{1}{p}\tilde{N}(v)\\
                       &=\left(\frac{1}{2}-\frac{1}{p}\right)a\|\nabla v \|_2^2 +\left(\frac{1}{2}-\frac{1}{p}\right)\int_{\mathbb{R}^N}(V(x)+\lambda)v^2dx +\left(\frac{1}{4}-\frac{1}{p}\right)bA^2\|\nabla v \|_2^2 \\
                       & \geq \left(\frac{1}{2}-\frac{1}{p}\right)a\|\nabla v \|_2^2 +\left(\frac{1}{2}-\frac{1}{p}\right)\int_{\mathbb{R}^N}\lambda v^2dx +\left(\frac{1}{4}-\frac{1}{p}\right)b\|\nabla v \|_2^4 \\
                       & \geq \left(\frac{1}{2}-\frac{1}{p}\right)at_*^2\|\nabla v \|_2^2 +\left(\frac{1}{2}-\frac{1}{p}\right)t_*^2\int_{\mathbb{R}^N}\lambda v^2dx +\left(\frac{1}{4}-\frac{1}{p}\right)bt_*^4\|\nabla v \|_2^4 \\
                       & = I_{\infty,\lambda}(t_*v)-\frac{1}{p}N_{\infty,\lambda}(t_*v)=I_{\infty,\lambda}(t_*v)\\
                       & \geq I_{\infty,\lambda}(u_{\beta})=m_{\beta}+\frac{1}{2}\lambda \beta^2.
\endaligned
\end{equation}
In addition, we have
$$  \aligned   (1+\|V \|_{\infty}) m_c \geq  m_c + \frac{1}{2}\| V \|_{\infty}c^2 > m_{V,c}& \geq J_{\lambda}(v)-\frac{\lambda}{2}\gamma^2 + m_{\alpha} \\
                                                             & = J_{\lambda}(v) +  m_{\alpha} + \frac{\lambda}{2}(\alpha^2-c^2) \\
                                                             & \geq m_{\beta} +  m_{\alpha} + \frac{\lambda}{2}(\alpha^2 + \beta^2 - c^2). \endaligned  $$
Since $m_{\beta}\geq m_{\alpha}>m_c$, we deduce $$ \alpha^2+\beta^2 < c^2. $$
Now, we estimate $\lambda$. By (\ref{lambda bound}), we have

$$ \lambda \leq   \left(2q+\frac{8\|V\|_{\infty}}{N(p-2)-4}\right)\frac{m_c}{c^2} + \frac{p(2-N)+2N}{N(p-2)-4}\frac{b}{2c^2}A^2.$$
By (\ref{A bound}) and (\ref{V12}), there exist a constant $K(a,p,N)>0$ depending on $a$, $N$ and $p$, such that
$$ \frac{p(2-N)+2N}{N(p-2)}\frac{b}{2}A^2 \leq 2 b K(a,p,N) m_c^2. $$
Therefore, we may write
$$ \lambda \leq  \left(2q+\frac{8\|V \|_{\infty}}{N(p-2)-4}\right)\frac{m_c}{c^2} + 2 b K(a,p,N)\frac{m_c}{c^2} m_c.  $$
Then, by (\ref{lambda bound}) and $\emph{Lemma 2.4}$, we have
 $$ \aligned (1+\|V \|_{\infty})m_c &> m_{V,c} \\
                  &\geq \left(\frac{c}{\beta}\right)^q m_c + \left(\frac{c}{\alpha}\right)^q m_c- \frac{\lambda}{2}(c^2-\alpha^2-\beta^2)\\
                  &\geq  \left(\frac{c}{\beta}\right)^q m_c + \left(\frac{c}{\alpha}\right)^q m_c-{\left( \left(q+\frac{4\|V\|_{\infty}}{N(p-2)-4}\right)+Kb m_c  \right)m_c} \frac{c^2-\alpha^2-\beta^2}{c^2},\\ \endaligned $$
that is
$$1 + q+\frac{N(p-2)\|V \|_{\infty}}{N(p-2)-4}+Kb m_c  \geq  \left(\frac{c}{\beta}\right)^q  + \biggl(\frac{c}{\alpha}\biggr)^q +{\left(\left(q+\frac{4\|V \|_{\infty}}{N(p-2)-4}\right)+Kb m_c \right)}\biggl(\frac{\alpha^2}{c^2}+\frac{\beta^2}{c^2}\biggr).$$
If $L>q$, then, by elementary arguments,
$$ \min_{\{x,y>0,x+y\leq 1\}}\{ x^{-q/2}+y^{-q/2}+L(x+y) \} =\left\{ \aligned &(2+q)\biggl(\frac{2L}{q}\biggr)^\frac{q}{2+q},\ \biggl(\frac{q}{2L}\biggr)^{2/(2+q)}\leq \frac{1}{2},\\
                                                                        & 2\cdot 2^{q/2}+L,\ \ \ \ \ \ \ \ \  \biggl(\frac{q}{2L}\biggr)^{2/(2+q)}\geq \frac{1}{2}. \endaligned \right. $$
Thus, we obtain
$$ (2+q)\biggl(\frac{2L}{q}\biggr)^\frac{q}{2+q} \geq (2+q)2^{q/2}\geq(2+q)\left( 1+\frac{q}{2}\ln2\right)>\frac{3}{2}q+2,$$
and
$$  2\cdot 2^{q/2}+L \geq 2+q\ln2+q>\frac{3}{2}q+2.   $$
Therefore,
$$1 + q+\frac{N(p-2)\|V \|_{\infty}}{N(p-2)-4}+Kb m_c  >  \frac{3}{2}q + 2.$$
However, we deduce
$$Kbm_c < \frac{2(p-2)}{N(p-2)-4}-\frac{N(p-2)}{N(p-2)-4}\|V \|_{\infty}$$
provided $bm_c <K^{-1} \left(\frac{2(p-2)}{N(p-2)-4}-\frac{N(p-2)}{N(p-2)-4}\|V \|_{\infty}\right):=L $ which is a contradiction. As a consequence, we have $v_n \rightarrow v$ in $H^1(\mathbb{R}^N)$ strongly, so that $\| v \|_2=c$.

\noindent\textbf{Proof of Theorem 1.2:}\
From \emph{Lemma 2.5} and (\ref{equ:Root equation}), we have
 $$\| \nabla Z_c \|_2 =\left[ \frac{4a}{N(p-2)} \right]^{\frac{4}{N(p-2)-4}} \| Q_p \|_{2}^{\frac{2(p-2)}{N(p-2)-4}}c^{\frac{4N-2p(N-2)}{4-N(p-2)}}, $$
$$I_{\infty}(Z_c)=\frac{N(p-2-4)}{2N(p-2)} \left[ \frac{4a}{N(p-2)} \right]^{\frac{4}{N(p-2)-4}} \| Q_p \|_{2}^{\frac{2(p-2)}{N(p-2)-4}}c^{\frac{4N-2p(N-2)}{4-N(p-2)}}. $$
Then there exist a constant $K_1(a,N,p,c,V)$ such that, $$bm_c \leq \frac{3}{2}b I_{\infty}(Z_c) < L \text{\ if\ } b <K_1(a,N,p,c,V). $$
From \emph{Lemma 2.2 (3)}, there exist a constant $K_2(a,N,p,b)$ such that, $$bm_c < L \text{\ if\ } c > K_2(a,N,p,b). $$
Combining (\ref{equ:lambda V}) and \emph{Lemma 4.6}, we obtain $(v,\lambda) \in H^1(\mathbb{R}^N)\times \mathbb{R}^+$ is a solution to $(K_{V,c})$. \qed

\noindent\textbf{Proof of Theorem 1.3:} Similar to \emph{Proposition 4.3} and \emph{Lemma 4.6}, we prove the conclusions under assumption $(V2)$.

Since $$C_n=\int_{\mathbb{R}^3}V(x)v_n dx \leq \|V \|_{\frac{3}{2}}\|v_n \|_{6}^2 \leq S^2\|V \|_{\frac{3}{2}} \|\nabla v_n\|_2^2 $$
and
$$D_n=\int_{\mathbb{R}^3}W(x)v_n(\nabla v_n \cdot x)dx \leq \| W\|_3 \|v_n\|_{6}^2 \|\nabla v_n \|_2^2 \leq S\|W \|_3 \|\nabla v_n\|_2^2,$$
we obtain
\begin{equation}\label{new bound}
  a(3p-10)A_n +\frac{b}{2}(3p-14)A_n^2 \leq 12(p-2)m_c + 4 \|W \|_3S A_n +6\|V \|_{\frac{3}{2}}S^2 A_n + o(1).
\end{equation}
Thus, $(v_n)_n$ is bounded in $H^1(\mathbb{R}^3)$ since $a(3p-10)-4 \|W \|_3S -6\|V \|_{\frac{3}{2}}S^2 >0$ provided (\ref{V22}).
We can argue as in the prove of \emph{Proposition 4.3}, and see that, up to a subsequence, $\lambda_n \rightarrow \lambda >0 $ provided
$$ (p-2)|D| < (6-p)m_c. $$

Using (\ref{new bound}), it is possible to see that
$$ \aligned  (p-2)|D|& \leq (p-2)\|W \|_3 A_n \\
                     & \leq \|W \|_3 \frac{12(p-2)^2m_c}{a(3p-10)-4 \|W \|_3S -6\|V \|_{\frac{3}{2}}S^2}
                     & \leq (6-p)m_c,
\endaligned $$
if (V2) holds. Finally, (\ref{V22}) yields
$$  c \leq S^2 \|V\|_{\frac{3}{2}} A \leq S^2 \| V \|_{\frac{3}{2}} \frac{12(p-2)m_c}{a(3p-10)-4 \|W \|_3S -6\|V \|_{\frac{3}{2}}S^2} \leq \frac{4}{3}m_c.$$
In the meanwhile, the estimate on $\lambda$ (\ref{lambda bound}) turns into
$$\aligned \lambda c^2 & \leq \left(2q +\frac{2(6-p)\nu}{3p-10}+\frac{4(p-2)}{3(3p-10)} \right)m_c +\frac{6-p}{3p-10}\frac{b}{2}A^2\\
                       & \leq \left(2q +\frac{2(6-p)\nu}{3p-10}+\frac{4(p-2)}{3(3p-10)} \right)m_c +\frac{6-p}{3p-10}\frac{b}{2}\left(\frac{12(p-2)m_c}{a(3p-10)-4 \|W \|_3S -6\|V \|_{\frac{3}{2}}S^2}\right)^2\\
                       & \leq \left(2q +\frac{2(6-p)\nu}{3p-10}+\frac{4(p-2)}{3(3p-10)} \right)m_c +\frac{6-p}{3p-10}\frac{b}{2}\left(\frac{12(p-2)}{a(3p-10)}\right)^2 m_c^2.\endaligned $$
Similarly, from \emph{Lemma 3.8}, we deduce
$$\aligned (1+\nu)m_c &> m_{V,c}  \\
                      &\geq \left(\frac{c}{\beta}\right)^q m_c + \left(\frac{c}{\alpha}\right)^q m_c- \frac{\lambda}{2}(c^2-\alpha^2-\beta^2)\\
                      &\geq  \left(\frac{c}{\beta}\right)^q m_c + \left(\frac{c}{\alpha}\right)^q m_c-{\left( \left(\frac{(6-p)\nu}{3p-10}+\frac{32-4p}{3(3p-10)}\right)+K^*b m_c  \right)m_c} \frac{c^2-\alpha^2-\beta^2}{c^2}\\ \endaligned $$
where $K^* = \frac{6-p}{4(3p-10)}\left(\frac{12(p-2)}{a(3p-10)}\right)^2$. Recall that $\nu < \frac{2}{3}$, arguing as in the proof of \emph{Lemma 4.6}, we can get a contradiction provided
$$ K^* b m_c <\frac{2(p-2)}{3p-10}\left(\frac{2}{3}-\nu\right). $$
Thus, $K_1^*(a,p,V,c)$ and $K_2^*(a,p,V,b)$ could be found and this completes the prove . \qed

\section{Existence of solutions under radial assumption}
\noindent\textbf{Case 1:} $V(x)$ satisfies $(V4)$ while $ 2 \leq N \leq 3$, $p \in [2+\frac{8}{N},2^*)$, and $c>c^*$.

\noindent\textbf{Proof of Theorem 1.5:}
In fact, adopting the argument of \emph{Section 4} to $\Psi_c$ and $S_c^r$, it is possible to prove the existence of a bounded sequence $v_n \in H_{rad}^1(\mathbb{R}^N)$ such that
$$  I(v_n)\rightarrow M_{V,c},\ \ \ \  \nabla_{S_c \cap H_{rad}^1(\mathbb{R}^N)}I(v_n)\rightarrow 0,  $$
and the Lagrange multipliers
$$ \lambda_n := -\frac{DI(v_n)[v_n]}{c^2}  $$
admits a subsequence converging to some $\lambda >0$, thanks to (\ref{V11})(\ref{V12}) or (\ref{V22}). We may assume that $\lim\limits_{n\rightarrow \infty} \| \nabla v_n \|_2^2 =A^2$ and $v_n \rightharpoonup v$ as $n \rightarrow \infty$. Then the weak limit $v$ solves
$$ -(a+bA^2)\triangle v +(V(x)+\lambda)v=|v|^{p-2}v \ \ \text{in}\ \mathbb{R}^N.   $$
Since the embedding $H_{rad}^1(\mathbb{R}^N) \hookrightarrow L^p(\mathbb{R}^N)$ is compact for $N\geq2$ and $p \in (2,2^*)$, we have $v_n \rightarrow v$ in $L^p(\mathbb{R}^N)$ strongly. As a consequence, using $\lambda>0$, we deduce
$$ v_n:=(-(a+bA^2)\triangle + V + \lambda )^{-1} \Bigg(|v_n|^{p-2}v_n +(\lambda - \lambda_n)v_n + b(A^2-\| \nabla v_n \|_2^2)\triangle v_n\Bigg) \rightarrow v $$
strongly in $H^1(\mathbb{R}^N)$, so $\|v\|_2=c$ and $(v, \lambda)$ is solution to $(K_{V,c})$.

\noindent\textbf{Case 2:} $V(x)$ satisfies $(V5)$ while $N\geq 2$, $2+\frac{4}{N} < p <\min\{2+\frac{8}{N},2^* \}$, and $c>c_1$.

\noindent\textbf{Proof of Theorem 1.6:}\ Similarly, we could get a sequence $v_n \in H_{rad}^1(\mathbb{R}^N)$ such that
$$  I(v_n)\rightarrow M_{V,c},\ \ \ \  \nabla_{S_c \cap H_{rad}^1(\mathbb{R}^N)}I(v_n)\rightarrow 0, \ \ P(v_n)\rightarrow 0. $$
Since $I$ is coercive on $S_c^r$ now, so that $v_n$ is bounded.
Without loss of generality, we may assume $$v_n \rightharpoonup v  \text{ in } H_{rad}^1(\mathbb{R}^N),\ v_n \rightarrow v \text{ in } L^p(\mathbb{R}^N) \text{ for } p \in (2,2^*),$$
and the Lagrange multipliers $\lambda_n$ satisfy
$$ \aligned - \lambda_n c^2 =& a\|\nabla v_n \|_2^2 + b\| \nabla v_n \|_2^4 + \int_{\mathbb{R}^N}V(x)v_n dx - \| v_n \|_p^p\\
                            =& \frac{N(p-2)-2p}{N(p-2)}(a\|\nabla v_n \|_2^2 + b\|\nabla v_n \|_2^4) + \int_{\mathbb{R}^N} \left(\frac{p}{N(p-2)} \langle \nabla V(x) \cdot x \rangle +V(x)\right)v_n^2dx  \endaligned $$
and admits a subsequence converging to some $\lambda \geq 0$.

If $v=0$, then
$$ \lim\limits_{n\rightarrow\infty}\int_{\mathbb{R}^N}\langle \nabla V(x)\cdot x  \rangle v_n^2 dx =0, $$
so, from $P(v_n) \rightarrow 0$ , we get $\| \nabla v_n \|_2^2 \rightarrow 0$ which contradicts $I(v_n) \rightarrow M_{V,c}>0$. Thus, $v \neq 0$ and $\lambda > 0$.
Meanwhile, if $\lim\limits_{n\rightarrow \infty} \| \nabla v_n \|_2^2 =A^2$, then $v$ solves
$$ -(a+bA^2)\triangle v +(V(x)+\lambda)v=|v|^{p-2}v \ \ \text{in}\ \mathbb{R}^N.   $$
As a consequence, using $\lambda>0$, we deduce
$$ v_n:=(-(a+bA^2)\triangle + V + \lambda )^{-1} \Bigg(|v_n|^{p-2}v_n +(\lambda - \lambda_n)v_n + b(A^2-\| \nabla v_n \|_2^2)\triangle v_n\Bigg) \rightarrow v $$
strongly in $H^1(\mathbb{R}^N)$, so $\|v\|_2=c$ and $(v, \lambda)$ is solution to $(K_{V,c})$.


\vskip8pt

\begin{thebibliography}{99}
\baselineskip=14pt
\bibitem{Arosio1997} A. Arosio, S. Panizzi: On the well-posedness of the Kirchhoff string. Tans. Am. Math. Soc. 348, 305-330 (1996).
\bibitem{Bartsch2021} T. Bartsch,  R. Molle, M. Rizzi, M. Verzini: Normalized solutions of mass supercritical Schr\"{o}dinger equations with potential. Comm. Partial Differ. Equ. 46 , no. 9, 1729-1756 (2021).
\bibitem{BartschJFA} T. Bartsch, N. Soave, A natural constraint approach to normalized solutions of nonlinear Schrödinger equations and systems, J. Funct. Anal. 272(12), 4998–5037 (2017).
\bibitem{Benci}V. Benci, G. Cerami: Positive solutions of some nonlinear elliptic problems in exterior domains. Arch. Rational Mech. Anal. 99, no. 4, 283–300 (1987).
\bibitem{Bernstein1940} S. Bernstein: Sur une class d'\'{e}quations fonctionnelles aux d\'{e}riv\'{e}es partielles. Bull. Acad. Sci. URSS. S\'{e}r. Math. 4, 17-26 (1940).
\bibitem{Ca2001} M. Cavalcanti, V. Cavalcanti, J. Soriano: Global existence and uniform decay rates for the Kirchhoff-Carrier equation with nonlinear dissipation. Adv. Differ. Equ. 6, 701-730 (2001).
\bibitem{CST2021} S.T. Chen, X.H. Tang: Normalized Solutions for Nonautonomous Schrödinger Equations on a Suitable Manifold. J. Geom. Anal. 30, 1637–1660 (2020).
\bibitem{ChenAMO} S.T. Chen, V.D. R\u{a}dulescu, X.H. Tang: Normalized Solutions of Nonautonomous Kirchhoff Equations: Sub- and Super-critical Cases. Appl. Math. Optim. 84, 773–806 (2021).
\bibitem{Spagnolo1992} P. D'Ancona, S. Spagnolo: Global solvability for the degenerate Kirchhoff equation with real analytic data. Invent. Math. 108, 247-262 (1992).
\bibitem{Ding2014} Y.H. Ding, X.X. Zhong: Normalized solution to the Schrödinger equation with potential and general nonlinear term: Mass super-critical case.
J. Diff. Equa. 334, 194-215 (2022).
\bibitem{Santos2014} G. Figueiredo, J.R. Santos: Multiplicity and concentration behavior of positive solutions for a Schr\"{o}dinger-Kirchhoff type problem via penalization method, Calc. Var. 20, 389-415 (2014).
\bibitem{Guo2015} Z.J. Guo: Ground states for Kirchhoff equations without compact condition. J. Differ. Equ. 259, 2884-2902 (2015).
\bibitem{Hu} T. Hu, C. Tang: Limiting behavior and local uniqueness of normalized solutions for mass critical Kirchhoff equations, Calc. Var. 60, 210 (2021).
\bibitem{Ikoma} N. Ikoma, Y. Miyamoto: Stable standing waves of nonlinear Schr¨odinger equations with potentials and general nonlinearities. Calc. Var. 59(2), 48, (2020).
\bibitem{Jeanjean1997} L. Jeanjean: Existence of solutions with prescribed norm for semilinear elliptic equations. Nonlinear Anal. 28(10), 1633–1659 (1997).
\bibitem{Kirchhoff1883} G. Kirchhoff: Mechanik, Teubner, Leipzig, (1883).
\bibitem{YeLi2019} G.B. Li, H.Y. Ye: On the concentration phenomenon of $L^2$-subcritical constrained minimizers for a class of Kirchhoff equations with potentials. J. Diff. Equa. 266, 7101-7123 (2019).
\bibitem{Li2021} G.B. Li, X. Luo, T. Yang: Normalized solutions to a class of Kirchhoff equations with Sobolev critical exponent.  arXiv:2103.08106v1 (2021).
\bibitem{Liyuhua2018} Y.H. Li, X.C. Hao, J.P. Shi: The existence of constrained minimizers for a class of nonlinear Kirchhoff–Schrödinger equations with doubly critical exponents in dimension four. Nonlinear Analysis. 186 99-112 (2019).
\bibitem{J.Lions1978} J.L. Lions: On some questions in boundary value problems of mathmatical physics. North-Holland Math. Stud. 30, 284-346 (1978).
\bibitem{Mao2009} A.M. Mao, Z.T. Zhang: Sign-changing and multiple solutions of Kirchhoff type problems without the P.S. condition, Nonlinear Anal. 70, 1275-1287 (2009).
\bibitem{Mao2022} A.M. Mao, S. Mo: Ground state solutions to a class of critical Schrödinger problem, Advances in Nonlinear Analysis vol. 11, no. 1, 96-127 (2022).
\bibitem{Molle2021} R. Molle, G. Riey, G. Verzini: Existence of normalized solutions to mass supercritical Schr\"{o}dinger equations with negative potential, arXiv:2104.12834v2. (2021)
\bibitem{Qi2022} S.J. Qi, W.M. Zou: Exact Number of Positive Solutions for the Kirchhoff Equation. SIAM Journal on Mathematical Analysis. vol.54, no.5, 5424-5446 (2022).
\bibitem{Shibata2014} M. Shibata: Stable standing waves of nonlinear Schr¨odinger equations with a general nonlinear term. Manuscripta Math. 143, 221-237, (2014).
\bibitem{Weinstein} M.I. Weinstein: Nonlinear Schr\"{o}dinger equations and sharp interpolation estimates. Commun. Math. Phys. 87, 567-576 (1983).

\bibitem{Xe2016} Q.L. Xie, S.W. Ma, X. Zhang: Bound state solutions of Kirchhoff type problems with critical exponent. J. Differ. Equ. 261(2), 890-924 (2016).
\bibitem{yang} Z. Yang: A new observation for the normalized solution of the Schr¨odinger equation. Arch. Math. 115 (3), 329-338 (2020).
\bibitem{Ye2014} H.Y. Ye: The sharp existence of constrained minimizers for a class of nonlinear Kirchhoff equations. Math. Methods
Appl. Sci. doi:10.1002/mma.3247 (2014)
\bibitem{Ye2015} H.Y. Ye: The existence of normalized solutions for $L^2$-critical constrained problems related to Kirchhoff equations. Z. Angew. Math. Phys. 66, 1483–1497 (2015).
\bibitem{Ye2016} H.Y. Ye: The mass concentration phenomenon for $L^2$-critical constrained problems related to Kirchhoff equations. Z. Angew. Math. Phys. 67, 29 (2016).
\bibitem{Zeng2017} X.Y. Zeng, Y.M. Zhang: Existence and uniqueness of normalized solutions for the Kirchhoff equation. Appl. Math. Lett. 74, 52-59 (2017).
\bibitem{Zeng2021} X.Y. Zeng, J.J. Zheng, Y.M. Zhang, X.X. Zhong: Positive normalized solution to the Kirchhoff equation with general nonlinearities. arXiv:2112.10293v1, (2021).
\bibitem{Zhang2022} P.H. Zhang, Z.Q. Han: Normalized ground states for Kirchhoff equations in $R^3$ with a critical nonlinearity. J. Math. Phys. 63, 021505 (2022).
\bibitem{Zhong2021} X.X. Zhong, W.M. Zou: A new deduce of the strict binding inequality and its application: Ground state normalized solution to Schr\"{o}dinger equations with potential. arXiv:2107.12558, (2021).
\bibitem{Zou2014} X.M. He, W.M. Zou: Ground states for nonlinear kirchhoff equations with critical growth. Annali di Matematica 193, 473-500 (2014).


\end{thebibliography}
\end{document}